# Degeneration of $l$-adic Eisenstein classes and of the elliptic polylog

Annette Huber huber@math.uni-muenster.de
Guido Kings kings@math.uni-muenster.de
Math. Institut, Einsteinstr. 62, 48 149 Münster, Germany

## Contents

# Introduction

Let $F = \mathbb{Q}(\mu_N)$ be a cyclotomic field. By Borel's deep theorem we know the ranks of all higher $K$-groups $K_*(F)$. It is

$$\dim K_n(F) \otimes \mathbb{Q} = \begin{cases} 0 & n > 0 \text{ even,} \\ r_2 & n > 1 \text{ odd.} \end{cases}$$

(We will not consider the more elementary and exceptional cases $n = 0, 1$.) However, Borel's proof does not produce explicit elements in $K_{2k+1}(F)$. Soulé nearly has achieved this. For each prime number $l$ and a primitive root of unity $\omega \in F$, he constructed a projective system of *cyclotomic elements* in $K_{2k+1}(F, \mathbb{Z}/l^r)$. It induces an element in $K_{2k+1}(F) \otimes \mathbb{Z}_l$ . The main consequence was the surjectivity of the higher Chern classes for $k \geq 1$

$$K_{2k+1}(F) \otimes \mathbb{Z}_l \to H^1(F, \mathbb{Z}_l(k+1))$$

where the right hand side is continuous Galois cohomology. It was conjectured that Soulé's elements for different $l$ are induced by the same element in $K_{2k+1}(F)$. This is indeed the case. More precisely, we give a new proof of the following main comparison theorem:

**Theorem (Beilinson, Deligne, Huber, Wildeshaus, [HuW]).** *For each root of unity $\omega \in F = \mathbb{Q}_l(\mu_N)$, ($N \geq 3$) with $\omega$ different from $1$ and for $k \geq 1$, there is an element in $K_{2k+1}(F)$ whose image in $K_{2k+1}(F) \otimes \mathbb{Z}_l$ and hence in $H^1(F, \mathbb{Z}_l(k+1))$ agrees with Soulé's cyclotomic elements for all $l$. Moreover the Beilinson regulator of these elements in Deligne cohomology can be computed and is given by higher logarithm functions.*

For the precise formulae we refer to theorem 1.4.4. These formulae were conjectured by Bloch and Kato ([BlK] 6.2). For $(l, N) = 1$ a weaker form was shown by Gros using syntomic cohomology. The general case was settled by constructing the motivic polylogarithm in [HuW]. The alternative proof which we are going to present in this paper is a lot less demanding on the technical side and quite short.

Our interest in the comparison theorem comes from its arithmetic applications. It is needed in attacking the conjectures on special values of $\zeta$-functions of more generally $L$-functions of abelian number fields. The Tamagawa number conjecture of Bloch and Kato was proved by themselves for the motives $\mathbb{Q}(n)$ with $n > 1$, i.e., the Riemann-$\zeta$-function, under the assumption of the above result ([BlK]). Another proved case is a version of the Lichtenbaum conjecture for abelian number fields ([KNF]).

We now want to sketch the line of the argument. The first key idea is to construct mixed motives for cyclotomic fields via modular curves. Anderson ($k = 0$) and Harder ($k \geq 0$) were the first to show that it is possible to construct mixed motives in this way. In fact, this paper is a sequel of [HuK] where we showed that Harder's construction has a $K$-theoretic analogue. So let $M$ be the modular curve over $F$ (see 1.1) and $j : M \to \overline{M}$ its compactification. Let $i : \operatorname{Spec} F \to \overline{M}$ the inclusion of the cusp $\infty$. Finally let $\mathcal{H}$ be the Tate-module of the universal elliptic curve over $M$. It is an $l$-adic sheaf on $M$. The *Eisenstein symbol* gives elements in

$$H^1_{\text{et}}(M, \operatorname{Sym}^k \mathcal{H}(1)) = H^1_{\text{et}}(M_{\overline{F}}, \operatorname{Sym}^k \mathcal{H}(1))^{G(\overline{F}/F)} .$$

Over the complex numbers they are related to Eisenstein series. The functor $i^* Rj_*$ - a derived version of degeneration at infinity - induces a little diagram

$$\begin{array}{c} H^1_{\text{et}}(M, \operatorname{Sym}^k \mathcal{H}(1)) \\ \downarrow \\ 0 \to H^1_{\text{et}}(F, i^* j_* \operatorname{Sym}^k \mathcal{H}(1)) \to H^1_{\text{et}}(F, i^* Rj_* \operatorname{Sym}^k \mathcal{H}(1)) \to H^0_{\text{et}}(F, i^* R^1 j_* \operatorname{Sym}^k \mathcal{H}(1)) \end{array}$$

In other words: an Eisenstein symbol whose residue in

$$H^0_{\text{et}}(\operatorname{Spec} F, i^* R^1 j_* \operatorname{Sym}^k \mathcal{H}(1)) = \mathbb{Q}_l$$

vanishes, induces an element of

$$H^1_{\text{et}}(\operatorname{Spec} F, i^* j_* \operatorname{Sym}^k \mathcal{H}(1)) = H^1(F, \mathbb{Q}_l(k+1)) ,$$

i.e., in the Galois cohomology group of $F$. We call these elements Harder-Anderson elements because this construction can be shown to be equivalent to Harder's original one. Now $K$-theory enters: Beilinson has constructed the Eisenstein symbol as element of $K$-groups of powers of the universal elliptic curve. The cup-product construction (see 1.2) in $l$-adic cohomology is a different way of writing down the above diagram. Phrased like this, the construction works in $K$-theory as well.

We now have to compute the degeneration of Eisenstein symbols with residue zero. In the Hodge theoretic analogue this can be done by direct computation (see [HuK] thm. 8.1). In the $l$-adic case, which is the subject of this paper, we use a different approach. We make use of the machinery of polylogarithms. All Eisenstein symbols appear as components of the *cohomological polylog* $\mathcal{P}\mathrm{ol}^{\mathrm{coh}}$ (see A.3.3) on the universal elliptic curve, more precisely as linear combinations of its fibres at torsion sections. If $\psi$ is a formal linear combination of torsion sections, we denote the corresponding linear combination of fibres by $\psi^* \mathcal{P}\mathrm{ol}^{\mathrm{coh}}$. If $\psi^* \mathcal{P}\mathrm{ol}^{\mathrm{coh}}$ has residue zero at infinity, the degeneration $i^* j_* \psi^* \mathcal{P}\mathrm{ol}^{\mathrm{coh}}$ is equal to the degeneration of the usual *elliptic polylog* $i^* j_* \psi^* \mathcal{P}\mathrm{ol}(-1)$. This degeneration is known and is equal to the *classical polylog*. Its value at torsion sections is again known and given by a linear combination of cyclotomic elements in $l$-adic cohomology. This finishes the computation. A different idea for the computation of the degeneration of $l$-adic Eisenstein classes was communicated to us by Harder.

The second subject of our paper (as is already apparent from the above sketch of our proof) is the systematic study of the degeneration of the elliptic into the classical polylogarithm induced by the degeneration of an elliptic curve into $\mathbb{G}_m$. Although the main results were known before, we found it helpful to introduce the concept of the cohomological polylogarithm in order to understand precisely what is going on. The two lines of thought, the connection of the cohomological polylogarithm and its degenerations with Eisenstein classes and the degeneration of the elliptic polylog into the classical polylog are given in the appendices C.1.1 and B.1.3, C.2.2. Unfortunately, we found it hard to quote the results from the literature. E.g., the comparison theorem C.2.2 is implicit in [BeL] and [W] but not stated. In a couple of cases (B.1.3 and C.1.1) it took us a while to understand the arguments in [BeL]. Given the nature of our application we had to be precise about factors and hence normalizations. In the end, we decided to make the paper as self-contained as possible. If it has doubled its length, it has hopefully also made it easier to read. Decisive was the urge to demonstrate that this is not at all a difficult theory.

**Overview:** Chapter 0 contains notations on sheaves and the like. In chapter 1 we first review Beilinson's Eisenstein symbol. Section 1.2 contains the main results of this article. The proof of the main theorem is given in chapter 2. Appendix A reviews the construction of the elliptic and the classical polylogarithm in unified form. In appendix B we study the degeneration of the elliptic into the classical polylog. The last appendix C treats the relation of the cohomological polylogarithm with the Eisenstein symbol. The appendices in the second part are meant to be (nearly) independent of the main text. Appendix B and C are independent of each other. Readers not familiar with polylogs should read section 0 and 1, then appendix A and then go on in the main text. We hope that even specialist may find the detailed proofs in appendix B and C useful.

**Acknowledgments:** We thank T. Scholl and J. Wildeshaus for discussions. The first author also wants to thank the DFG for financial support during the time the article was prepared.

# 0 Notations and conventions

Let $l$ be a prime. All schemes are separated and of finite type over some regular scheme of dimension zero or one and $l$ is invertible on them.

By *$l$-adic sheaves* we mean constructible $\mathbb{Q}_l$-sheaves. The category of $l$-adic sheaves on a scheme $X$ is denoted $\mathrm{Sh}(X)$. We abbreviate $\mathrm{Hom}_{\mathrm{Sh}(X)} = \mathrm{Hom}_X$. By $\underline{\mathrm{Hom}}_X$ we denote the internal Hom with values in $\mathrm{Sh}(X)$. In particular, if $\mathcal{H}$ is a lisse $l$-adic sheaf, then $\mathcal{H}^\vee = \underline{\mathrm{Hom}}_X(\mathcal{H}, \mathbb{Q}_l)$ is its dual. By the standard map $\mathbb{Q}_l \to \mathcal{H}^\vee \otimes \mathcal{H}$ we mean the dual of the evaluation map $\mathcal{H} \otimes \mathcal{H}^\vee \to \mathbb{Q}_l$. Whenever we have a map of $l$-adic sheaves $\mathcal{H} \otimes \mathcal{F} \to \mathcal{G}$, then we mean by the induced map the composition

$$\mathcal{F} \to \mathcal{H}^\vee \otimes \mathcal{H} \otimes \mathcal{F} \to \mathcal{H}^\vee \otimes \mathcal{G} \ .$$

We will have to use infinite direct products and projective limits of $l$-adic sheaves. These are to be understood formally. We really mean projective systems of sheaves.

Let $D^*(X)$ be the "derived" category of $l$-adic sheaves in the sense of Ekedahl [E]. It is a triangulated category with a canonical $t$-structure whose heart is $\mathrm{Sh}(X)$. Note that $D(X)$ is *not* the derived category, i.e., morphisms do not compute Ext-groups. However, for $\mathcal{F}, \mathcal{G} \in \mathrm{Sh}(X)$, we still have

$$\mathrm{Hom}_X(\mathcal{F}, \mathcal{G}) = \mathrm{Hom}_{D(X)}(\mathcal{F}, \mathcal{G}) \ , \ \mathrm{Ext}^1_X(\mathcal{F}, \mathcal{G}) = \mathrm{Hom}_{D(X)}(\mathcal{F}, \mathcal{G}[1]) \ .$$

For $\mathcal{F} \in \mathrm{Sh}(F)$, we put

$$H^i_{\mathrm{et}}(X, \mathcal{F}) = \mathrm{Hom}_{D(X)}(\mathbb{Q}_l, \mathcal{F}[i]) \ .$$

By [Hu] 4.1, it agrees with continuous étale cohomology as defined by Jannsen [J]. If $f : X \to Y$ is a morphism, then there are functors

$$Rf_*, Rf_! : D^+(X) \to D^+(Y)$$
$$Rf^*, Rf^! : D^+(Y) \to D^+(X)$$

For their properties we refer to [E] theorem 6.3.

**Definition 0.0.1.** *Suppose $X$ is a smooth variety over a field and $U$ an open subvariety, such that the closed complement $Y$ is also smooth and of pure codimension* 1*.*

**a)** *The* residue map

$$\mathrm{res}: H^q_{\mathrm{et}}(U, \mathbb{Q}_l(k)) \to H^{q-1}_{\mathrm{et}}(Y, \mathbb{Q}_l(k-1))$$

*is defined as the connecting morphism of the long exact cohomology sequence for the triangle $i_*\mathbb{Q}_l(-1)[-2] \to \mathbb{Q}_l \to Rj_*\mathbb{Q}_l$ (note that $i^!\mathbb{Q}_l \cong i^*\mathbb{Q}_l(-1)[-2]$ by purity).*

**b)** *For sheaves $\mathcal{F}$ and $\mathcal{G}$ on $U$ and a morphism in the derived category $f: \mathcal{F} \to \mathcal{G}[q]$, we define the* residue of $f$ *by*

$$\mathrm{res}(f): i^*j_*\mathcal{F} \to i^*Rj_*\mathcal{F} \to i^*Rj_*\mathcal{G}[q] \to i^*R^1j_*\mathcal{G}[q-1]$$

*(note that by assumption $Rj_*$ has cohomological dimension 1).*

It is easy to see that these definitions agree when $\mathcal{F} = \mathbb{Q}_l$, $\mathcal{G} = \mathbb{Q}_l(k)$. We will also use motivic cohomology.

**Definition 0.0.2.** *If $X$ is a regular scheme, then we call*

$$H^i_{\mathcal{M}}(X, j) := \mathrm{Gr}^\gamma_j K_{2j-i}(X) \otimes \mathbb{Q}$$

*the* motivic cohomology *of $X$.*

**Remark:** For schemes of finite type over a number field we can also use the category of horizontal sheaves and its "derived" category as defined in [Hu]. In this case, $K_*(X)$ has to be replaced by the direct limit of the $K_*(\mathfrak{X})$ where $\mathfrak{X}$ runs through all $\mathbb{Z}$-smooth models of $X$. The arguments in the article work without any changes.

# 1 The Eisenstein symbol and the cup product construction

In this section we recall the Eisenstein symbol constructed by Beilinson. We then explain our main construction. It yields elements of $H^1_{\mathcal{M}}(\mathrm{Spec}\,\mathbb{Q}(\mu_N), k+1)$. Finally we formulate the main theorem - the computation of the regulators of these elements - and its consequences.

## 1.1 A residue sequence in $K$-theory

Let $M$ be the modular curve parameterizing elliptic curves with full level-$N$-structure and let $N \geq 3$. It is a geometrically connected variety defined over $B := \operatorname{Spec} \mathbb{Q}(\mu_N)$. We also fix a primitive root of unity $\zeta \in \mathbb{Q}(\mu_N)$. Let $E$ be the universal elliptic curve with level-$N$-structure over $M$. Let $\overline{M}$ be the compactification of $M$ and $\widetilde{E}$ the Néron model of $E$ over $\overline{M}$. Let $e : \overline{M} \to \widetilde{E}$ be the unit section. Its connected component in $\widetilde{E}$ is denoted by $\widetilde{E}^0$. Let $\mathrm{Cusp} = \overline{M} \smallsetminus M$ be the subscheme of cusps. The standard-$N$-gon over $B$ with level-$N$-structure $\mathbb{Z}/N \times \mathbb{Z}/N \to \mathbb{G}_m \times \mathbb{Z}/N$ via $(a, b) \mapsto (\zeta^a, b)$ *defines* a section $\infty : B \to \mathrm{Cusp}$. We have a diagram

$$\begin{array}{ccccc} E & \xrightarrow{j} & \widetilde{E} & \longleftarrow & \widetilde{E}_{\mathrm{Cusp}} \\ {\scriptstyle\pi}\downarrow & & {\scriptstyle\widetilde{\pi}}\downarrow & & {\scriptstyle\widetilde{\pi}}\downarrow \\ M & \xrightarrow{j'} & \overline{M} & \longleftarrow & \mathrm{Cusp} \end{array} .$$

The Eisenstein symbol will give elements in the motivic cohomology of

$$E^k = E \times_M \cdots \times_M E$$

the $k$-fold relative fibre product. The inversion $\iota$ operates on $E$ and the symmetric group $\mathfrak{S}_k$ operates on $E^k$. This induces an operation of the semi-direct product $\mu_2^k \rtimes \mathfrak{S}_k$ on $E^k$. Denote by $\varepsilon$ the character

$$\varepsilon : \mu_2^k \rtimes \mathfrak{S}_k \to \mu_2$$

which is the multiplication on $\mu_2^k$ and the sign-character on $\mathfrak{S}_k$. The localization sequence for the pair $((\widetilde{E}^0)^k, E^k)$ gives

$$H^{k+1}_{\mathcal{M}}(E^k, k+1)(\varepsilon) \to H^k_{\mathcal{M}}((\widetilde{E}^0)^k_{\mathrm{Cusp}}, k)(\varepsilon),$$

where $(\dots)(\varepsilon)$ is the $\varepsilon$ eigenspace.

Recall that $\widetilde{E}_{\mathrm{Cusp}}$ is non-canonically isomorphic to $\mathbb{G}_m \times \mathbb{Z}/N \times \mathrm{Cusp}$. This would induce an isomorphism (via residue at zero) $H^k_{\mathcal{M}}((\widetilde{E}^0)^k_{\mathrm{Cusp}}, k)(\varepsilon) \cong H^0_{\mathcal{M}}(\mathrm{Cusp}, 0)$. However, this isomorphism is *not* equivariant under $-\mathrm{id} \in \mathrm{Gl}_2(\mathbb{Z}/N)$. In order to make the isomorphism canonical, we introduce the $\mu_2$-torsor on Cusp given by

$$\mathrm{Isom} = \mathrm{Isom}(\mathbb{G}_m, \widetilde{E}^0_{\mathrm{Cusp}}) .$$

Over $\infty$, we have $\mathrm{Isom}_\infty = \mu_{2,\infty}$ by the very choice of $\infty$.

The isomorphism (given by residue at zero on $\mathbb{G}_m$)

$$H^k_{\mathcal{M}}((\widetilde{E}^0)^k_{\mathrm{Isom}}, k)(\varepsilon) \cong H^0_{\mathcal{M}}(\mathrm{Isom}, 0)$$

*is* $\mathrm{Gl}_2$-equivariant.

Composing with the pull-back to $(\widetilde{E}^0)^k_{\mathrm{Isom}}$ we get a map

$$\mathrm{res}^k : H^{k+1}_{\mathcal{M}}(E^k, k+1)(\varepsilon) \to H^0_{\mathcal{M}}(\mathrm{Isom}, 0)$$

which we call the *residue map.* It takes values in the subspace $H^0_{\mathcal{M}}(\mathrm{Isom}, 0)^{(k)}$ on which $-1 \in \mu_2$ acts via $(-1)^k$.

Explicitly, we have (using $\infty$)

$$\mathrm{Isom}(\overline{\mathbb{Q}}) = P(\mathbb{Z}/N) \setminus \mathrm{Gl}_2(\mathbb{Z}/N)$$

where $P = \begin{pmatrix} * & * \\ 0 & 1 \end{pmatrix} \subset \mathrm{Gl}_2$ with trivial action of $\mathrm{Gal}(\overline{\mathbb{Q}}/\mathbb{Q}(\mu_N))$. In this correspondence, the identity matrix corresponds to $\mathrm{id} \in \mu_{2,\infty} = \mathrm{Isom}_\infty$.

**Definition 1.1.1.** *Let*

$$\mathbb{Q}[\mathrm{Isom}]^{(k)} = \{f : \mathrm{Gl}_2(\mathbb{Z}/N) \to \mathbb{Q} \mid f(ug) = f(g) \text{ for } u \in P(\mathbb{Z}/N) \text{ and } f(-\mathrm{id}\, g) = (-1)^k f(g)\}\ ,$$

*the space of formal linear combination of points of* Isom *on which* $\mu_2$ *operates by* $(-1)^k$*. We can identify*

$$H^0_{\mathcal{M}}(\mathrm{Isom}, 0)^{(k)} = \mathbb{Q}[\mathrm{Isom}]^{(k)}\ .$$

**Remark:** By choice of a section, $\mathbb{Q}[\mathrm{Isom}]^{(k)}$ is isomorphic to $\mathbb{Q}[\mathrm{Cusp}]$ but again this isomorphism cannot be made $\mathrm{Gl}_2$-equivariant and will not be used.

## 1.2 The Eisenstein symbol

**Proposition 1.2.1 (Beilinson).** *For* $k \geq 0$*, there is a map, called the* Eisenstein symbol

$$\mathrm{Eis}^k : \mathbb{Q}[\mathrm{Isom}]^{(k)} \to H^{k+1}_{\mathcal{M}}(E^k, k+1)(\varepsilon)$$

*It is a splitting of the residue map, i.e.,*

$$\mathrm{res}^k \circ \mathrm{Eis}^k = \mathrm{id}\ .$$

*Proof.* This is [Be] Theorem 7.3 or [Sch-Sch] section 7. □

**Remark:** Beilinson's construction really works for all elliptic curves over regular base schemes.

**Definition 1.2.2.** *Let* $\mathrm{Eis}_l^k$ *be the composition of* $\mathrm{Eis}^k$ *with the regulator map*

$$r_l : H^{k+1}_{\mathcal{M}}(E^k, k+1)(\varepsilon) \to H^{k+1}_{\mathrm{et}}(E^k, \mathbb{Q}_l(k+1))(\varepsilon) \ .$$

Let $B_k(x)$ be the $k$-th Bernoulli polynomial, defined by

$$\frac{te^{tx}}{e^t - 1} = \sum_{k=0}^{\infty} B_k(x) \frac{t^k}{k!} \ .$$

For $\overline{x} \in \mathbb{R}/\mathbb{Z}$ let $B_k(\overline{x}) = B_k(x)$ where $x$ is the representative of $\overline{x}$ in $[0,1)$.

**Definition 1.2.3.** *The* horospherical map

$$\varrho^k : \mathbb{Q}[E[N]] \longrightarrow \mathbb{Q}[\mathrm{Isom}]^{(k)}$$

*maps* $\psi : E[N] \to \mathbb{Q}$ *to*

$$\varrho^k(\psi)(g) = \frac{N^k}{k!(k+2)} \sum_{\substack{t=(t^{(1)},t^{(2)}) \\ \in (\mathbb{Z}/N)^2}} \psi(g^{-1}t) B_{k+2}\left(\frac{t^{(2)}}{N}\right) .$$

This map is well-defined, $\mathrm{Gl}_2$-equivariant and surjective. It has a kernel. In particular:

**Lemma 1.2.4.** *Let* $\mathbb{Q}[E[N] \smallsetminus 0]$ *be the* $\mathbb{Q}$*-vector space of maps* $\psi : E[N] \smallsetminus 0 \to \mathbb{Q}$ *of degree* $0$*. Then*

$$\varrho : \mathbb{Q}[E[N] \smallsetminus 0] \to \mathbb{Q}[\mathrm{Isom}]^{(k)}$$

*is still surjective.*

*Proof.* Use the $\mathrm{Gl}_2$-translates of the element $\phi_\infty \in \mathbb{Q}[E[N] \smallsetminus 0]$ written out in [HuK] lemma 7.6. □

**Remark:** Beilinson first constructs the composition $\mathrm{Eis}^k \circ \varrho^k$ and then deduces the existence of the section.

## 1.3 The cup-product construction

We now construct mixed Dirichlet motives starting from the Eisenstein symbol. This construction already appears in [HuK] Section 4, where the relation to Harder-Anderson motives is also clarified.

Let $\infty = \operatorname{Spec}\mathbb{Q}(\mu_N) \subset \mathrm{Cusp}$ be as above.

**Definition 1.3.1.** *Let* $\mathrm{res}^k_\infty$ *be the composition of* $\mathrm{res}^k$ *with the restriction to* $\infty$*. Let*

$\mathbb{Q}[\mathrm{Isom}\smallsetminus\infty]^{(k)} = \left\{ f \in \mathbb{Q}[\mathrm{Isom}]^{(k)} \text{(cf. 1.1.1)} \mid f(g) = 0 \text{ for } g \in \pm P(\mathbb{Z}/N) \right\}$ .

*Note that* $\varrho^k\psi \in \mathbb{Q}[\mathrm{Isom}\smallsetminus\infty]^{(k)}$ *is equivalent to* $\mathrm{res}^k_\infty \mathrm{Eis}^k(\varrho^k\psi) = 0$*. Let*

$$\phi_\infty \in \mathbb{Q}[\mathrm{Isom}]^{(k)}$$

*be given by*

$$\phi_\infty(g) = \begin{cases} 1 & g \in P(\mathbb{Z}/N), \\ (-1)^k & g \in -P(\mathbb{Z}/N), \\ 0 & \text{else.} \end{cases}$$

To define the space of Harder Anderson elements, consider:

$$H^{k+1}_{\mathcal{M}}(E^k, k+1) \xrightarrow{\cup \mathrm{Eis}^k(\phi_\infty)} H^{2k+2}_{\mathcal{M}}(E^k, 2k+2) \xrightarrow{\pi_*} H^2_{\mathcal{M}}(M, k+2) \xrightarrow{\mathrm{res}^k_\infty} H^1_{\mathcal{M}}(\infty, k+1) \ .$$

**Definition 1.3.2.** *For* $k \geq 1$ *let*

$$\mathrm{Dir} : \mathbb{Q}[\mathrm{Isom}\smallsetminus\infty]^{(k)} \to H^1_{\mathcal{M}}(\infty, k+1)$$
$$\phi \mapsto \mathrm{res}^k_\infty \circ \pi_*(\mathrm{Eis}^k(\phi) \cup \mathrm{Eis}^k(\phi_\infty))$$

$\mathrm{Dir}\left(\mathbb{Q}[\mathrm{Isom}\smallsetminus\infty]^{(k)}\right)$ *is called the space of* Harder-Anderson elements. *Let*

$$\mathrm{Dir}_l : \mathbb{Q}[\mathrm{Isom}\smallsetminus\infty]^{(k)} \to H^1_{\mathrm{et}}(\infty, \mathbb{Q}_l(k+1))$$

*be the composition of* Dir *with the* $l$*-adic regulator.*

**Remark:** We call elements in $H^1_{\mathcal{M}}(\operatorname{Spec}\mathbb{Q}(\mu_N), k+1)$ Dirichlet motives because their Hodge regulator is related to Dirichlet series. Dir attaches a Dirichlet motive to a linear combination of cusps. As shown in [HuK] 5.2 and 6.4., the cup-product construction, or more precisely its image in $l$-adic or Hodge cohomology can be translated into Harder's construction in [Ha] 4.2, whence the name.

## 1.4 The main theorem

Our aim is to compute the image of $\mathrm{Dir}(\phi)$ for $\phi \in \mathbb{Q}[\mathrm{Isom} \smallsetminus \infty]^{(k)}$ under the $l$-adic regulator. Recall that

$$\mathrm{Dir}_l(\phi) \in H^1_{\mathrm{et}}(\infty, \mathbb{Q}_l(k+1)) = H^1_{\mathrm{et}}(\mathrm{Spec}\, \mathbb{Q}(\mu_N), \mathbb{Q}_l(k+1))$$

We will use the following explicit representation of this group.

**Lemma 1.4.1.** *Recall that $\zeta$ is a fixed primitive $N$-th root of unity. There is an isomorphism*

$$\left( \varprojlim_{r \geq 1} H^1_{\mathrm{Gal}}(\mathbb{Q}(\mu_N), \mathbb{Z}/l^r(k+1)) \right) \otimes \mathbb{Q}_l \cong H^1_{\mathrm{et}}(\mathrm{Spec}\, \mathbb{Q}(\mu_N), \mathbb{Q}_l(k+1)) \ .$$

*Proof.* The choice of $\zeta$ induces an isomorphism between étale and Galois cohomology. All $\mathbb{Z}/l^r(k+1)$ are finite, hence $\lim^1$ is zero. □

Soulé has constructed elements in $K$-theory with coefficients ([Sou]) whose image in $l$-adic cohomology is the following:

**Definition 1.4.2.** *For $u \in \mathbb{Z}/N, u \neq 0$ we call the class*

$$c^k(\zeta^u) := \left( \sum_{\alpha^{l^r} = \zeta^u} [1 - \alpha] \otimes (\alpha^N)^{\otimes k} \right)_r \in H^1_{\mathrm{et}}(\mathrm{Spec}\, \mathbb{Q}(\mu_N), \mathbb{Q}_l(k+1))$$

*(under the identification of 1.4.1) Soulé-Deligne or cyclotomic element.*

After these preparations we can formulate our main theorem. Recall that the horospherical map (see definition 1.2.3)

$$\varrho^k : \mathbb{Q}[E[N] \smallsetminus 0] \to \mathbb{Q}[\mathrm{Isom}]^{(k)}$$

is surjective.

**Theorem 1.4.3 (The main theorem).** *Let $k \geq 1$ and*

$$\psi \in \mathbb{Q}[E[N] \smallsetminus 0] \ .$$

*Assume that $\varrho^k \psi \in \mathbb{Q}[\mathrm{Isom} \smallsetminus \infty]^{(k)}$, i.e., that*

$$\mathrm{res}^k_\infty \mathrm{Eis}^k(\varrho^k \psi) = 0 \ .$$

*Then*

$$\mathrm{Dir}(\varrho^k \psi) \in H^1_{\mathcal{M}}(\mathrm{Spec}\, \mathbb{Q}(\mu_N), k+1)$$

*is defined (cf. definition 1.3.2).*

**a)** *The $l$-adic regulator of* $\mathrm{Dir}(\varrho^k\psi)$ *in* $H^1_{\mathrm{et}}(\operatorname{Spec}\mathbb{Q}(\mu_N), \mathbb{Q}_l(k+1))$ *is given by*

$$\mathrm{Dir}_l(\varrho^k\psi) = \frac{(-1)^{k+1}}{k!N}\sum_{t\in\mathbb{Z}/N}\psi(t,0)c^k(\zeta^t) \ .$$

**b)** *The Hodge regulator of* $\mathrm{Dir}(\varrho^k\psi)$ *in* $H^1_{\mathcal{H}}(\operatorname{Spec}\mathbb{Q}(\mu_N)\otimes_{\mathbb{Q}}\mathbb{C}, \mathbb{R}(k+1))$ *is given by*

$$r_{\mathcal{H}}\,\mathrm{Dir}(\varrho^k\psi) = (-1)^{k+1}N^{k-1}\left(\sum_{t\in\mathbb{Z}/N}\psi(t,0)Li_{k+1}(\sigma\zeta^t)\right)_{\sigma:\mathbb{Q}(\mu_N)\hookrightarrow\mathbb{C}}$$

*via the identification*

$$H^1_{\mathcal{H}}(\operatorname{Spec}\mathbb{Q}(\mu_N)\otimes_{\mathbb{Q}}\mathbb{C}, \mathbb{R}(k+1)) \cong \prod_{\sigma:\mathbb{Q}(\mu_N)\hookrightarrow\mathbb{C}}\mathbb{C}/\mathbb{R}(k) \ ,$$

*and where* $Li_{k+1}(x) := \sum_{n>0}\frac{x^n}{n^{k+1}}$.

This theorem will be proved in section 2.4. Part b) is theorem 8.1 in [HuK]. All normalizations agree.

**Remark:** Note that the simple shape of the above formulae is due to the parameterization of $\mathbb{Q}[\mathrm{Isom}]^{(k)}$ via $\varrho^k$. For $\psi\in\mathrm{Ker}(\varrho^k)$, the formulae imply relations between values of higher logarithm functions respectively between Soulé-Deligne elements.

As an immediate consequence we get:

**Corollary 1.4.4 (Bloch-Kato conjecture 6.2).** *For* $u\neq 0$ *consider the element*

$$\psi^k_u := \frac{(-1)^{k+1}}{N^{k-1}}(u,0) - \frac{(-1)^{k+1}N^2}{1-N^{k+1}}\sum_{v\neq 0}(u,v)\in\mathbb{Q}[E[N]\smallsetminus 0].$$

*By definition of the horospherical map (1.2.3) and the distribution relation for Bernoulli numbers,* $\varrho^k\psi^k_u\in\mathbb{Q}[\mathrm{Isom}\smallsetminus\infty]^{(k)}$. *Hence*

$$\mathrm{Dir}(\varrho^k\psi^k_u)\in H^1_{\mathcal{M}}(\operatorname{Spec}\mathbb{Q}(\mu_N), k+1)$$

*is defined. Its regulators are given by*

$$r_{\mathcal{H}}\,\mathrm{Dir}(\varrho^k\psi^k_u) = [Li_{k+1}(\sigma(\zeta^u))]_{\sigma:\mathbb{Q}(\mu_N)\to\mathbb{C}}$$

*and*

$$\mathrm{Dir}_l(\varrho^k \psi_u^k) = \frac{1}{N^k k!} c^k(\zeta^u) \ .$$

**Remark:** The sign differs from the one given in [HuW] because of a different normalization in the identification of $H^1_{\mathcal{H}}(B \otimes_{\mathbb{Q}} \mathbb{C}, \mathbb{R}(k+1))$ with a product of $\mathbb{C}/\mathbb{R}(k)$'s in [W] Part II Thm 3.6 p. 222.

The existence of elements of $H^1_{\mathcal{M}}(\mathrm{Spec}\,\mathbb{Q}(\mu_N), k+1)$ with Hodge-regulator as above has been proved by Beilinson [Be]. The compatibility of his elements with Soulé's was a conjecture of Bloch and Kato ([BlK] Conjecture 6.2). A local version has previously been proved in the case $(l, N) = 1$ by Gros ([G] Thm 2.4) using the syntomic regulator. The general case was settled by Huber and Wildeshaus ([HuW]) using a completely different method due to Beilinson and Deligne.

The implications are of course the same as in [HuW]: Completion of the proof of Tamagawa number conjecture of Bloch and Kato for the motives $\mathbb{Q}(n+1)$, $n \geq 1$ given in [BlK] 6.1 ii; a version of the Lichtenbaum conjecture on special values of $\zeta$-functions for abelian number fields ([KNF] Thm 6.4); the proof that Soulé's elements in $K$-theory with coefficients are induced by elements in $K$-theory itself. For precise formulations we refer to [HuW] 9.6 to 9.8. Moreover, the above theorem gives a direct computation of the $l$-adic version of Harder's motives in [Ha] 4.2. For the precise relation we refer to [HuK] 6.4.

In the course of the proof of the main theorem 1.4.3 a), we will also show the following assertion, which is interesting in itself:

**Theorem 1.4.5.** *We use the notation of 2.1. Let* $k \geq 1$. *Let* $\phi \in \mathbb{Q}[\mathrm{Isom} \smallsetminus \infty]^{(k)}$, *i.e.,* $\mathrm{Eis}_l^k(\phi)$ *in the kernel of* $\mathrm{res}^k_\infty$. *Let* $\mathcal{E}^k_\phi$ *be the sheaf representing* $\mathrm{Eis}^k_l(\phi) \in \mathrm{Ext}^1_M(\mathbb{Q}_l, \mathrm{Sym}^k \mathcal{H}(1))$. *Finally let* $\psi \in \mathbb{Q}[E[N] \smallsetminus 0]$ *with* $\varrho^k\psi = \phi$. *Then*

$$i'^* j'_* \mathcal{E}^k_\phi \in \mathrm{Ext}^1_\infty(\mathbb{Q}_l, \mathbb{Q}_l(k+1))$$

*and it is given by* $\frac{(-1)^{k+1}}{k! N} \sum_{t \in \mathbb{Z}/N} \psi(t, 0) c^k(\zeta^t)$ *(cf. 1.4.2).*

The proof will also be given in section 2.3.

## 2 Proof of the main theorem

The strategy is to link the Eisenstein-symbol to the elliptic polylog and the cyclotomic elements to the cyclotomic polylog. Then the degeneration of the elliptic polylog into the cyclotomic polylog is used to show the theorem.

### 2.1 Translation to a degeneration theorem

We want to reduce the computation of $\mathrm{Dir}_l(\phi)$ to the study of the degeneration of Eisenstein classes at $\infty$, see theorem 2.1.4.

**Definition 2.1.1.** *Define*

$$\mathcal{H} := (R^1\pi_*\mathbb{Q}_l)^\vee = \underline{\mathrm{Hom}}(R^1\pi_*\mathbb{Q}_l, \mathbb{Q}_l)$$

*to be the dual of the relative first cohomology group on $M$. Hence, in every fibre, $\mathcal{H}$ is the Tate module.*

Before we formulate the theorem, we need to understand $\mathcal{H}$ in a small neighbourhood of $\infty$. Let $S$ be the completion of $\overline{M}$ at $\infty$, $\eta$ its generic point. The closed point is $\infty = B$. We have the diagram

$$\begin{array}{ccccc} E_\eta & \xrightarrow{j} & \widetilde{E} & \xleftarrow{i} & \mathbb{G}_m \times \mathbb{Z}/N \\ {\scriptstyle \pi_\eta}\downarrow & & {\scriptstyle \widetilde{\pi}}\downarrow & & \downarrow{\scriptstyle \widetilde{\pi}} \\ \eta & \xrightarrow{j'} & S & \xleftarrow{i'} & B \end{array}$$

**Lemma 2.1.2.** *The sheaf $\mathcal{H}_\eta$ (pull-back of $\mathcal{H}$ from $M$ to $\eta$) has a canonical filtration*

$$0 \to \mathbb{Q}_l(1) \xrightarrow{\partial} \mathcal{H}_\eta \xrightarrow{p} \mathbb{Q}_l \to 0\ .$$

*$i'^*j'_*\mathcal{H}_\eta$ is canonically isomorphic to the Tate module of the special fibre $\widetilde{E}^0_s = \mathbb{G}_m$ which is identified with $\mathbb{Q}_l(1)$ via residue at zero. $\partial$ is normalized such that $i'^*j'_*(\partial) = \mathrm{id}$. The map $p$ is normalized such that the intersection pairing $\mathcal{H}_\eta \otimes \mathcal{H}_\eta \to \mathbb{Q}_l(1)$ induces $p(1)$ on $\mathbb{Q}_l(1) \otimes \mathcal{H}_\eta \to \mathbb{Q}_l(1)$.*

The proof of this (well-known) lemma is given in B.1.1.

**Remark:** In the usual complex parameterization, let $\{1, \tau\}$ the standard basis for $H_1(E_\tau)$, where $E_\tau = \mathbb{C}/(\mathbb{Z} + \tau\mathbb{Z})$. Then $\partial(1) = 1$, $p(\tau) = 1$.

**Corollary 2.1.3.** $\mathrm{Sym}^k(\partial)$ *induces an isomorphism*

$$\mathbb{Q}_l(n) \xrightarrow{\sim} i'^* j'_* \mathrm{Sym}^n \mathcal{H}_\eta \ .$$

*The projections via* $p$

$$\mathrm{Sym}^n \mathcal{H}_\eta \to \mathrm{Sym}^{n-1} \mathcal{H}_\eta$$

*induce isomorphisms*

$$i'^* R^1 j'_* \mathrm{Sym}^n \mathcal{H}_\eta \cong i'^* R^1 j'_* \mathrm{Sym}^{n-1} \mathcal{H}_\eta \cong \cdots \cong i'^* R^1 j'_* \mathbb{Q}_l \cong \mathbb{Q}_l(-1) \ .$$

*Proof.* This is shown by induction on $n$ by considering the long exact sequence attached to

$$0 \to \mathbb{Q}_l(n) \to \mathrm{Sym}^n \mathcal{H}_\eta \to \mathrm{Sym}^{n-1} \mathcal{H}_\eta \to 0 \ .$$

□

We now reformulate the cup-product construction in sheaf theoretic terms. Via the Leray spectral sequence for $R\pi_*$ and Poincaré duality for $\mathcal{H}$ we get

$$H^{k+1}_{\mathrm{et}}(E^k, \mathbb{Q}_l(k+1))(\varepsilon) = H^1_{\mathrm{et}}(M, \mathrm{Sym}^k \mathcal{H}(1)) \ .$$

**Theorem 2.1.4.** *For* $\phi \in \mathbb{Q}[\mathrm{Isom} \smallsetminus \infty]^{(k)}$*, let* $\mathcal{E}^k_\phi$ *be the sheaf representing* $\mathrm{Eis}^k_l(\phi) \in H^1_{\mathrm{et}}(M, \mathrm{Sym}^k \mathcal{H}(1)) = \mathrm{Ext}^1_M(\mathbb{Q}_l, \mathrm{Sym}^k \mathcal{H}(1))$*. Let* $i', j'$ *be as above. Then the sequence*

$$0 \to i'^* j'_* \mathrm{Sym}^k \mathcal{H}(1) \to i'^* j'_* \mathcal{E}^k_\phi \to \mathbb{Q}_l \to 0$$

*is exact and via the isomorphism (2.1.3)*

$$i'^* j'_* \mathrm{Sym}^k \mathcal{H}(1) \xrightarrow{\sim} \mathbb{Q}_l(k+1) \ ;$$

*it represents*

$$\mathrm{Dir}_l(\phi) = r_l \, \mathrm{res}^k_\infty(\mathrm{Eis}^k(\phi) \cup \mathrm{Eis}^k(\phi_\infty)) \in H^1_{\mathrm{et}}(\infty, \mathbb{Q}_l(k+1)).$$

We need a couple of lemmas:

**Lemma 2.1.5.** *We have a commutative diagram*

$$\begin{array}{ccc} H^{2k+2}_{\mathrm{et}}(E^k \times E^k, \mathbb{Q}_l(2k+2))(\varepsilon,\varepsilon) & \xrightarrow{\Delta^*} & H^{2k+2}_{\mathrm{et}}(E^k, \mathbb{Q}_l(2k+2)) \\ \cong\downarrow & & \downarrow\pi_* \\ H^2_{\mathrm{et}}(M, \mathrm{Sym}^k \mathcal{H}(1) \otimes \mathrm{Sym}^k \mathcal{H}(1)) & \longrightarrow & H^2_{\mathrm{et}}(M, \mathbb{Q}_l(k+2)) \end{array}$$

*where the last map is induced by the intersection pairing.*

*Proof.* The intersection pairing

$$\mathcal{H} \otimes \mathcal{H} \to \mathbb{Q}_l(1)$$

translates via Poincaré duality $\mathcal{H}^\vee \cong \mathcal{H}(-1)$ into into the cup-product pairing

$$\mathcal{H}^\vee \otimes \mathcal{H}^\vee \to R^2\pi_*\mathbb{Q}_l \xrightarrow{\sim} \mathbb{Q}_l(-1)\ .$$

Note that it uses the same trace map which defines $\pi_*$ on cohomology. □

**Lemma 2.1.6.** *Let $f : \mathbb{Q}_l \to \mathcal{F}[1]$ and $g : \mathcal{F} \to \mathcal{G}[1]$ be morphisms in $D(M)$ where $\mathcal{F}$ and $\mathcal{G}$ are sheaves. We assume that the map $i'^*j'_*\mathbb{Q}_l \to i'^*R^1j'_*\mathcal{F}$ induced by $f$ vanishes. Then the residue of $g \circ f : \mathbb{Q}_l \to \mathcal{G}[2]$ (cf. definition 0.0.1) is given by the composition*

$$\mathbb{Q}_l \to (i'^*j'_*\mathcal{F})[1] \to (i'^*R^1j'_*\mathcal{G})[1].$$

*In other words by the push-out of the short exact sequence*

$$0 \to i'^*j'_*\mathcal{F} \to i'^*j'_*\mathcal{E}_f \to \mathbb{Q}_l \to 0$$

*induced by $f$ via the map $i'^*j'_*\mathcal{F} \to i'^*R^1j'_*\mathcal{G}$ induced by $g$.*

*Proof.* Recall how the residue of $g \circ f$ is defined: it is the composition

$$\mathbb{Q}_l \to i'^*Rj'_*\mathbb{Q}_l \xrightarrow{i'^*Rj'_*f} i'^*Rj'_*\mathcal{F}[1] \xrightarrow{i'^*Rj'_*g} i'^*Rj'_*\mathcal{G}[2] \to i'^*R^1j'_*\mathcal{G}[1]\ .$$

By assumption $\mathbb{Q}_l \to i'^*Rj'_*\mathcal{F}[1]$ factors through $\mathbb{Q}_l \to i'^*j'_*\mathcal{F}[1]$. □

*Proof.* (of theorem 2.1.4) We have

$$\mathrm{Dir}_l(\phi) := r_l(\mathrm{res}^k_\infty(\mathrm{Eis}^k(\phi) \cup \mathrm{Eis}^k(\phi_\infty))) = \mathrm{res}^k_\infty(\mathrm{Eis}^k_l(\phi) \cup \mathrm{Eis}^k_l(\phi_\infty))\ .$$

We write the cup-product as composition

$$\begin{array}{ccccc} \mathbb{Q}_l & \xrightarrow{\mathrm{Eis}^k_l(\phi)\otimes\mathrm{id}} & \mathrm{Sym}^k\,\mathcal{H}(1)[1] \otimes \mathbb{Q}_l & \xrightarrow{\mathrm{id}\otimes\mathrm{Eis}^k_l(\phi_\infty)} & \mathrm{Sym}^k\,\mathcal{H}(1) \otimes \mathrm{Sym}^k\,\mathcal{H}(1)[2] \\ & & & & \downarrow \\ & & & & \mathbb{Q}_l(k+2)[2] \end{array}$$

Now we apply the previous lemma with $\mathcal{F} = \mathrm{Sym}^k\,\mathcal{H}(1)$ and $\mathcal{G} = \mathbb{Q}_l(k+2)$. It remains to show that the map induced by $\mathrm{Eis}^k_l(\phi_\infty)$

$$\mathbb{Q}_l(k+1) \xrightarrow{\sim} i'^*j'_*\,\mathrm{Sym}^k\,\mathcal{H}(1) \to i'^*R^1j'_*\mathbb{Q}_l(k+2) \xrightarrow{\sim} \mathbb{Q}_l(k+1)$$

is the identity. We have a commutative diagram

$$\begin{array}{ccc}
i'^*j'_*\operatorname{Sym}^k\mathcal{H}(1)\otimes i'^*j'_*\mathbb{Q}_l & \xrightarrow{\mathrm{id}\otimes\mathrm{res}^k_\infty(\mathrm{Eis}^k_l(\phi_\infty))} & i'^*j'_*\operatorname{Sym}^k\mathcal{H}(1)\otimes i'^*R^1j'_*\operatorname{Sym}^k\mathcal{H}(1) \\
{=}\downarrow & & \downarrow \\
i'^*j'_*(\operatorname{Sym}^k\mathcal{H}(1)\otimes\mathbb{Q}_l) & \longrightarrow & i'^*R^1j'_*(\operatorname{Sym}^k\mathcal{H}(1)\otimes\operatorname{Sym}^k\mathcal{H}(1)) \\
 & & \downarrow \\
 & & i'^*R^1j'_*\mathbb{Q}_l(k+2)
\end{array}$$

The two vertical arrows are given by the cup-product. The top-most map induces the identity by assumption on $\phi_\infty$ (see definition 1.3.1). That the composition of the two vertical maps on the right is also the identity follows from our normalization in definition 2.1.2 with respect to the intersection pairing. □

## 2.2 Connection with the cohomological polylog

We now want to compare the Eisenstein classes with the elliptic cohomological polylog. The theory is developed in [BeL], [W]. For a detailed exposition we refer to the appendix.

Let $\mathcal{L}\mathrm{og}_{E/M}$ on $E$ and $\mathcal{L}\mathrm{og}_{\mathbb{G}_m/B}$ on $\mathbb{G}_m$ be the elliptic and the classical logarithmic (pro)-sheaf respectively. Furthermore let

$$\mathcal{P}\mathrm{ol}_{\mathbb{G}_m/B}\in\operatorname{Ext}^1_{\mathbb{G}_m\smallsetminus e(B)}(\widetilde{\pi}^*\mathbb{Q}_l(1),\mathcal{L}\mathrm{og}_{\mathbb{G}_m/B})\ ,$$
$$\mathcal{P}\mathrm{ol}_{E/M}\in\operatorname{Ext}^1_{E\smallsetminus e(M)}(\pi^*\mathcal{H},\mathcal{L}\mathrm{og}_{E/M})$$

be the classical and the elliptic polylog respectively. We will also need a variant of this.

**Definition 2.2.1.** *The* cohomological polylogarithm $\mathcal{P}\mathrm{ol}^{\mathrm{coh}}_{E/M}$ *is the class*

$$\mathcal{P}\mathrm{ol}^{\mathrm{coh}}_{E/M}\in\operatorname{Ext}^1_{E\smallsetminus e(M)}(\mathbb{Q}_l,\pi^*\mathcal{H}^\vee\otimes\mathcal{L}\mathrm{og}_{E/M}(1))$$

*obtained by tensoring* $\mathcal{P}\mathrm{ol}_{E/M}$ *with* $\mathcal{H}^\vee$ *and pull-back with the standard map* $\mathbb{Q}_l\to\mathcal{H}^\vee\otimes\mathcal{H}$*. The representing sheaf will also be denoted by* $\mathcal{P}\mathrm{ol}^{\mathrm{coh}}_{E/M}$*.*

Recall (e.g. A.2.6) that we have isomorphisms for all torsion sections $t$

$$\varrho_t : t^* \mathcal{L}\mathrm{og}_{E/M} \cong e^* \mathcal{L}\mathrm{og}_{E/M} = \prod_{k \geq 0} \mathrm{Sym}^k \mathcal{H} \ ,$$

$$\varrho_t : t^* \mathcal{L}\mathrm{og}_{\mathbb{G}_m/B} \cong e^* \mathcal{L}\mathrm{og}_{\mathbb{G}_m/B} = \prod_{k \geq 0} \mathbb{Q}_l(k) \ .$$

For convenience, we make the following definition.

**Definition 2.2.2. a)** *For $\psi = \sum q_i t_i \in \mathbb{Q}[E[N] \smallsetminus 0]$ put*

$$\psi^* \mathcal{P}\mathrm{ol}^{\mathrm{coh}}_{E/M} = \sum q_i \varrho_{t_i} (t_i^* \mathcal{P}\mathrm{ol}^{\mathrm{coh}}_{E/M})$$

*in $\mathrm{Ext}^1_M(\mathbb{Q}_l, \mathcal{H}^\vee \otimes e^* \mathcal{L}\mathrm{og}_{E/M}(1))$. Similarly we define $\psi^* \mathcal{P}\mathrm{ol}_{E/M}$.*

**b)** *Via the isomorphism*

$$\mathrm{Ext}^1_M(\mathbb{Q}_l, \mathcal{H}^\vee \otimes e^* \mathcal{L}\mathrm{og}_{E/M}(1)) = \prod_{k \geq 0} \mathrm{Ext}^1_M(\mathbb{Q}_l, \mathcal{H}^\vee \otimes \mathrm{Sym}^k \mathcal{H}(1))$$

*we denote by*

$$\left(\psi^* \mathcal{P}\mathrm{ol}^{\mathrm{coh}}_{E/M}\right)^k \in \mathrm{Ext}^1_M(\mathbb{Q}_l, \mathcal{H}^\vee \otimes \mathrm{Sym}^k \mathcal{H}(1))$$

*the $k$-th component of $\psi^* \mathcal{P}\mathrm{ol}^{\mathrm{coh}}_{E/M}$. Similarly we define $\left(\psi^* \mathcal{P}\mathrm{ol}_{E/M}\right)^k$.*

**Lemma 2.2.3.** *The operation of $\mathcal{H}$ on $\mathcal{L}\mathrm{og}_{E/M}$ (see appendix A.2) induces*

$$\mu^\vee : \mathcal{L}\mathrm{og}_{E/M} \to \pi^* \mathcal{H}^\vee \otimes \mathcal{L}\mathrm{og}_{E/M} \ .$$

*$e^* \mu^\vee$ has a section*

$$\mathrm{pr} : \mathcal{H}^\vee \otimes e^* \mathcal{L}\mathrm{og}_{E/M} \to e^* \mathcal{L}\mathrm{og}_{E/M}$$

*which on $\mathcal{H}^\vee \otimes \mathrm{Sym}^k \mathcal{H}$ is given by*

$$h^\vee \otimes h_0 \otimes \cdots \otimes h_k \mapsto \frac{1}{k+2} \sum_{j=0}^{k} h^\vee(h_j) h_0 \otimes \cdots \hat{h}_j \cdots h_k \ .$$

*It induces an isomorphism*

$$\mathrm{Ext}^1_M(\mathbb{Q}_l, \mathcal{H}^\vee \otimes e^* \mathcal{L}\mathrm{og}_{E/M}(1)) \to \mathrm{Ext}^1_M(\mathbb{Q}_l, e^* \mathcal{L}\mathrm{og}_{E/M}(1))$$

*It maps the $k+1$-st component on the left to the $k$-th component on the right.*

*Proof.* See [W], 3.19 b). Use the Leray spectral sequence and the known weights of the cohomology of $\mathrm{Sym}^k \mathcal{H}$ over $\overline{\mathbb{Q}}$. □

We can now formulate the theorem that makes the polylog so important in our context:

**Theorem 2.2.4 (Elliptic polylog at torsion sections).**
*Let $\psi \in \mathbb{Q}[E[N] \smallsetminus 0]$ and $k > 0$, then*

$$\mathrm{Eis}_l^k(\varrho^k \psi) = -N^{k-1} \,\mathrm{pr}\left(\psi^* \,\mathcal{P}\mathrm{ol}_{E/M}^{\mathrm{coh}}\right)^{k+1}$$

*as elements of*

$$\mathrm{Ext}^1_M(\mathbb{Q}_l, \mathrm{Sym}^k \mathcal{H}(1))$$

*with the projection* pr *in lemma 2.2.3.*

This is implicit in [W] p. 310, and [BeL] 2.1-2.2. For the detailed exposition see the appendix, C.1.1.

**Remark:** This theorem is also true for $k = 0$.

## 2.3 Degeneration of Eis and $\mathcal{P}\mathrm{ol}_{E/M}^{\mathrm{coh}}$

The aim of this section is to show:

**Theorem 2.3.1.** *Let $\psi \in \mathbb{Q}[E[N] \smallsetminus 0]$ with $\varrho^k\psi \in \mathbb{Q}[\mathrm{Isom} \smallsetminus \infty]^{(k)}$, i.e., such that*

$$\mathrm{Eis}_l^k(\varrho^k \psi) \in \mathrm{Ext}^1_M(\mathbb{Q}_l, \mathrm{Sym}^k \mathcal{H}(1))$$

*has residue 0 at $\infty$. Recall that our fixed choice of $\zeta$ induces an inclusion $\mathbb{Z}/N \subset \mathbb{G}_m(B)$ via $t \mapsto \zeta^t$. Then in $\mathrm{Ext}^1_B(\mathbb{Q}_l, \mathbb{Q}_l(k+1))$*

$$i'^* j'_* \,\mathrm{Eis}_l^k(\varrho^k \psi) = -N^{k-1} \sum_{t \in \mathbb{Z}/N} \psi(t, 0) \left(\varrho_t t^* \,\mathcal{P}\mathrm{ol}_{\mathbb{G}_m/B}(-1)\right)^{k+1} \ .$$

We first relate $i'^* j'_* \psi^* \,\mathcal{P}\mathrm{ol}_{E/M}^{\mathrm{coh}}$ to $i'^* j'_* \psi^* \,\mathcal{P}\mathrm{ol}_{E/M}$.

It is enough to restrict to the local situation around $\infty$, with notation as in 2.1. Let $\partial^\vee : \mathcal{H}_\eta^\vee \to \mathbb{Q}_l(-1)$ be the dual of the map $\partial$ defined in 2.1.2.

**Lemma 2.3.2.** *Let $\psi \in \mathbb{Q}[E[N] \smallsetminus 0]$ such that*

$$\operatorname{res}^k_\infty \left(\psi^* \mathcal{P}\mathrm{ol}^{\mathrm{coh}}_{E/\eta}\right)^{k+1} = 0 \ .$$

*Then the sequence*

$$0 \to i'^* j'_* (\mathcal{H}^\vee_\eta \otimes \mathrm{Sym}^{k+1} \mathcal{H}(1)) \to i'^* j'_* \left(\psi^* \mathcal{P}\mathrm{ol}^{\mathrm{coh}}_{E/\eta}\right)^{k+1} \to \mathbb{Q}_l \to 0$$

*is exact. Push-out by $i'^* j'_*(\partial^\vee \otimes \mathrm{id})$ is the extension*

$$0 \to i'^* j'_* \mathrm{Sym}^{k+1} \mathcal{H} \to i'^* j'_* \left(\psi^* \mathcal{P}\mathrm{ol}_{E/\eta}(-1)\right)^{k+1} \to \mathbb{Q}_l \to 0.$$

*Proof.* By definition of the cohomological polylog we have a commutative diagram (we forget about the index $E/\eta$ and $\eta$)

$$\begin{array}{ccccccccc}
0 & \longrightarrow & \pi^*\mathcal{H}^\vee \otimes \mathcal{L}\mathrm{og}(1) & \longrightarrow & \mathcal{P}\mathrm{ol}^{\mathrm{coh}} & \longrightarrow & \mathbb{Q}_l & \longrightarrow & 0 \\
 & & \downarrow{=} & & \downarrow & & \downarrow & & \\
0 & \longrightarrow & \pi^*\mathcal{H}^\vee \otimes \mathcal{L}\mathrm{og}(1) & \longrightarrow & \mathcal{H}^\vee \otimes \mathcal{P}\mathrm{ol} & \longrightarrow & \mathcal{H}^\vee \otimes \mathcal{H} & \longrightarrow & 0
\end{array}$$

We apply $\psi^*$ to it and take the $k+1$-st component. We get a commutative diagram

$$\begin{array}{ccccccccc}
0 & \longrightarrow & \mathcal{H}^\vee \otimes \mathrm{Sym}^{k+1}\mathcal{H}(1) & \longrightarrow & \left(\psi^* \mathcal{P}\mathrm{ol}^{\mathrm{coh}}\right)^{k+1} & \longrightarrow & \mathbb{Q}_l & \longrightarrow & 0 \\
 & & \downarrow{=} & & \downarrow & & \downarrow & & \\
0 & \longrightarrow & \mathcal{H}^\vee \otimes \mathrm{Sym}^{k+1}\mathcal{H}(1) & \longrightarrow & \mathcal{H}^\vee \otimes (\psi^* \mathcal{P}\mathrm{ol})^{k+1} & \longrightarrow & \mathcal{H}^\vee \otimes \mathcal{H} & \longrightarrow & 0 \\
 & & \downarrow{\partial^\vee} & & \downarrow{\partial^\vee} & & \downarrow{\partial^\vee} & & \\
0 & \longrightarrow & \mathrm{Sym}^{k+1}\mathcal{H} & \longrightarrow & (\psi^* \mathcal{P}\mathrm{ol}(-1))^{k+1} & \longrightarrow & \mathcal{H}(-1) & \longrightarrow & 0
\end{array}$$

We apply the functor $i'^* j'_*$ to this diagram. The first line remains exact by assumption on the residue. The last line remains exact because $\mathcal{P}\mathrm{ol}$ and $\psi^* \mathcal{P}\mathrm{ol}$ have residue zero by B.1.4. The composition of the right vertical maps is the identity. □

*Proof.* (of theorem 2.3.1) Recall that by lemma 2.2.3 the map pr induces an isomorphism

$$\mathrm{Ext}^1_M(\mathbb{Q}_l, \mathcal{H}^\vee \otimes e^* \mathcal{L}\mathrm{og}_{E/M}(1)) \xrightarrow{\mathrm{pr}} \mathrm{Ext}^1_M(\mathbb{Q}_l, e^* \mathcal{L}\mathrm{og}_{E/M}(1)) \ .$$

Hence for $\left(\psi^* \mathcal{P}\mathrm{ol}^{\mathrm{coh}}_{E/M}\right)^{k+1}$ and its image under pr to have residue zero is equivalent. Application of $i'^* j'_*$ to the formula in theorem 2.2.4 gives

$$i'^* j'_* \mathcal{E}^k_{E/M}(\varrho^k \psi) = -N^{k-1} i'^* j'_*(\mathrm{pr}) i'^* j'_* \psi^* \mathcal{P}\mathrm{ol}^{\mathrm{coh}}_{E/M} \ .$$

It is easy to check that both $i'^* j'_*(\mathrm{pr})$ and $i'^* j'_*(\partial^\vee \otimes \mathrm{id})$ are sections of

$$i'^* j'_*(\mu^\vee) : \mathbb{Q}_l(k+1) = i'^* j'_* \left(\mathrm{Sym}^k \mathcal{H}_\eta(1)\right) \to i'^* j'_* \left(\mathcal{H}^\vee_\eta \otimes \mathrm{Sym}^{k+1} \mathcal{H}_\eta(1)\right) \ .$$

But such a section is unique for weight reasons and hence

$$i'^* j'_*(\mathrm{pr}) = i'^* j'_*(\partial^\vee \otimes \mathrm{id}) \ .$$

Hence we get

$$\begin{aligned} i'^* j'_* \mathcal{E}^k_{E/M}(\varrho^k \psi) &= -N^{k-1} i'^* j'_*(\partial^\vee \otimes \mathrm{id}) i'^* j'_* \left(\psi^* \mathcal{P}\mathrm{ol}^{\mathrm{coh}}_{E/\eta}\right)^{k+1} \\ &\overset{2.3.2}{=} -N^{k-1} i'^* j'_* \left(\psi^* \mathcal{P}\mathrm{ol}_{E/\eta}(-1)\right)^{k+1} \\ &\overset{B.2.1}{=} -N^{k-1} \left(\psi^* i^* j_* \mathcal{P}\mathrm{ol}_{E/\eta}(-1)\right)^{k+1} \\ &= -N^{k-1} \left(\sum_{t \in \mathbb{Z}/N} \psi(t,0) \varrho_t t^* \mathcal{P}\mathrm{ol}_{\mathbb{G}_m/B}(-1)\right)^{k+1} \ . \end{aligned}$$

The last equality follows from the known degeneration of the elliptic polylog (see B.1.4): if $t$ is a torsion section of $\widetilde{E}$ meeting $\mathbb{G}_m \times \{v\}$, then

$$\begin{aligned} &t^* i^* j_* \mathcal{P}\mathrm{ol}_{E/\eta} = i'^* j'_* t^* \mathcal{P}\mathrm{ol}_{E/\eta} \cong t^* \mathcal{P}\mathrm{ol}_{\mathbb{G}_m/B} \ \text{ for } v = 0, \\ &t^* i^* j_* \mathcal{P}\mathrm{ol}_{E/\eta} = i'^* j'_* t^* \mathcal{P}\mathrm{ol}_{E/\eta} \ \text{ is split for } v \neq 0. \end{aligned}$$

□

## 2.4 End of proof

*Proof.* (of theorem 1.4.3 a) and theorem 1.4.5) We have shown

$$\mathrm{Dir}_l(\varrho^k \psi) \overset{2.1.4}{=} i'^* j'_* \mathcal{E}^k_{\varrho^k \psi} \overset{2.3.1}{=} -N^{k-1} \left(\sum_{t \in \mathbb{Z}/N} \psi(t,0) \varrho_t t^* \mathcal{P}\mathrm{ol}_{\mathbb{G}_m/B}(-1)\right)^{k+1} \ .$$

If $t \in \mathbb{Z}/N \setminus 0$, then the section of $\mathbb{G}_m$ associated to $t$ is $\zeta^t$, and we have

$$\left(\varrho_t t^* \mathcal{P}\mathrm{ol}_{\mathbb{G}_m/B}(-1)\right)^{k+1} = \frac{(-1)^k}{N^k k!} c^k(\zeta^t)$$
$$= \frac{(-1)^k}{N^k k!} \left( \sum_{\alpha^{l^r} = \zeta^{t(1)}} [1-\alpha] \otimes (\alpha^N)^{\otimes k} \right)_r .$$

This is [W] part II, chapter 4. Note that our $N$ is $d$ in loc. cit. whereas $N$ in loc. cit. is 1 in our case. □

**Remark:** There are two possible strategies for the Hodge theoretic counterpart. Everything we have said about lisse $l$-adic sheaves works also for admissible variations of $\mathbb{Q}$-Hodge structure. There is a Hodge theoretic version of $\mathcal{H}$, $\mathcal{L}$og and $\mathcal{P}$ol and the theorems remain valid. The only new ingredient is the explicit computation of the cyclotomic polylog at a torsion section. For this see [W] Part II Theorem 3.11. Hence the arguments of this chapter give a valid proof of b). On the other hand, there is a direct computational proof of b), see [HuK] Theorem 8.1.

# A Elliptic and classical polylogarithm

In this appendix we want to give a short review of the classical and the elliptic polylogarithm. We assemble the facts which are needed in the main text. Other important aspects are omitted. For example, we restrict to the $l$-adic setting. None of the material is new and it is only included to make our presentation as self-contained as possible. The classical polylogarithm was constructed by Deligne in [Del2]. The elliptic polylogarithm is due to Beilinson and Levin in [BeL] and we follow their presentation to a large extent. Another presentation of the material in this section can be found in Wildeshaus [W] IV, V.

## A.1 The logarithmic sheaf

Let $\pi : G \to S$ be an elliptic curve or a torus of relative dimension one and $e : S \to G$ be the unit section. In this section we introduce and study the logarithm sheaf.

**Definition A.1.1.** *Define*

$$\mathcal{H} := (R^1\pi_*\mathbb{Q}_l)^\vee = \underline{\mathrm{Hom}}(R^1\pi_*\mathbb{Q}_l, \mathbb{Q}_l)$$

*to be the dual of the relative first cohomology group on $S$. Hence, in every fibre, $\mathcal{H}$ is the Tate module.*

**Lemma A.1.2.** *For every lisse sheaf $\mathcal{F}$ on $S$, there is a canonical isomorphism*

$$R^q\pi_*\pi^*\mathcal{F} \cong (R^{2-q}\pi_!\mathbb{Q}_l(1))^\vee \otimes \mathcal{F} \cong (R^q\pi_*\mathbb{Q}_l) \otimes \mathcal{F} \ .$$

*Proof.* This is just the projection formula together with Poincaré duality:

$$\begin{aligned} R^q\pi_*\pi^*\mathcal{F} &\cong (R^{2-q}\pi_!\pi^*(\mathcal{F}^\vee))^\vee(-1) \\ &\cong (R^{2-q}\pi_!\mathbb{Q}_l(1) \otimes (\mathcal{F}^\vee))^\vee \cong (R^{2-q}\pi_!\mathbb{Q}_l(1))^\vee \otimes \mathcal{F}. \end{aligned}$$

□

Consider the exact sequence coming from the Leray spectral sequence for $R\pi_*$. With the above lemma it reads

$$0 \to \mathrm{Ext}^1_S(\mathbb{Q}_l, \mathcal{H}) \to \mathrm{Ext}^1_G(\mathbb{Q}_l, \pi^*\mathcal{H}) \to \mathrm{Hom}_S(\mathbb{Q}_l, \mathcal{H}^\vee \otimes \mathcal{H}) \to 0 \ .$$

This sequence is split by $e^*$.

**Definition A.1.3.** *Let $\mathcal{L}\mathrm{og}^{(1)}_{G/S}$ be the unique extension*

$$0 \to \pi^*\mathcal{H} \to \mathcal{L}\mathrm{og}^{(1)}_{G/S} \to \mathbb{Q}_l \to 0$$

*such that its image in $\mathrm{Hom}_S(\mathbb{Q}_l, \mathcal{H}^\vee \otimes \mathcal{H})$ is the standard morphism*

$$\mathbb{Q}_l \to \mathcal{H}^\vee \otimes \mathcal{H}$$

*and $e^*\mathcal{L}\mathrm{og}^{(1)}_{G/S}$ is split (the splitting is unique for weight reasons). Let*

$$\mathcal{L}\mathrm{og}^{(n)}_{G/S} := \mathrm{Sym}^n \mathcal{L}\mathrm{og}^{(1)}_{G/S}$$

*and*

$$\mathcal{L}\mathrm{og}_{G/S} := \varprojlim \mathcal{L}\mathrm{og}^{(n)}_{G/S}$$

*where the limit is taken with respect to the transition maps*

$$\mathrm{Sym}^{k+1} \mathcal{L}\mathrm{og}^{(1)}_{G/S} \to \mathrm{Sym}^{k+1}(\mathcal{L}\mathrm{og}^{(1)}_{G/S} \oplus \mathbb{Q}_l) \to \mathrm{Sym}^k \mathcal{L}\mathrm{og}^{(1)}_{G/S} \otimes \mathbb{Q}_l$$

*(the first map is induced by the identity and* $\mathcal{L}\mathrm{og}^{(1)}_{G/S} \to \mathbb{Q}_l$ *and the second map is the canonical projection in the symmetric algebra of a direct sum). Hence, there are exact sequences*

$$0 \to \pi^* \operatorname{Sym}^{n+1} \mathcal{H} \to \mathcal{L}\mathrm{og}^{(n+1)}_{G/S} \to \mathcal{L}\mathrm{og}^{(n)}_{G/S} \to 0.$$

*The (pro-) sheaf* $\mathcal{L}\mathrm{og}_{G/S}$ *is called the* logarithm sheaf. *The splitting of* $e^* \mathcal{L}\mathrm{og}^{(1)}$ *induces*

$$e^* \mathcal{L}\mathrm{og}_{G/S} \cong \prod_{k \geq 0} \operatorname{Sym}^k \mathcal{H} \ .$$

Very important for everything that follows, is the computation of the higher direct images of the logarithm sheaf. Consider the exact sequence

$$0 \to \pi^* \operatorname{Sym}^{k} \mathcal{H} \to \mathcal{L}\mathrm{og}^{(k)}_{G/S} \to \mathcal{L}\mathrm{og}^{(k-1)}_{G/S} \to 0.$$

**Lemma A.1.4. a)** *Let $G/S$ be an elliptic curve. Then*

$$\begin{array}{ll} R\pi_* = R\pi_! \ , & \pi_*\mathcal{L}\mathrm{og}^{(k)} = \operatorname{Sym}^k \mathcal{H} \ , \\ R^1\pi_*\mathcal{L}\mathrm{og}^{(k)} = \operatorname{Sym}^{k+1} \mathcal{H}(-1) \ , & R^2\pi_*\mathcal{L}\mathrm{og}^{(k)} = \mathbb{Q}_l(-1) \ . \end{array}$$

**b)** *Let $G/S$ be a torus. Then*

$$\begin{array}{ll} \pi_*\mathcal{L}\mathrm{og}^{(k)} = \operatorname{Sym}^k \mathcal{H} \ , & \pi_!\mathcal{L}\mathrm{og}^{(k)} = 0 \ , \\ R^1\pi_*\mathcal{L}\mathrm{og}^{(k)} = \mathbb{Q}_l(-1) \ , & R^1\pi_!\mathcal{L}\mathrm{og}^{(k)} = \operatorname{Sym}^{k+1} \mathcal{H}(-1) \ , \\ R^2\pi_*\mathcal{L}\mathrm{og}^{(k)} = 0 \ , & R^2\pi_!\mathcal{L}\mathrm{og}^{(k)} = \mathbb{Q}_l(-1) \ . \end{array}$$

**c)** *In both cases, the boundary maps*

$$\begin{array}{r} \pi_*\mathcal{L}\mathrm{og}^{(k-1)} \to R^1\pi_* \operatorname{Sym}^k \mathcal{H} = \mathcal{H}^\vee \otimes \operatorname{Sym}^k \mathcal{H} \ , \\ R^1\pi_!\mathcal{L}\mathrm{og}^{(k-1)} \to R^2\pi_! \operatorname{Sym}^k \mathcal{H} \end{array}$$

*are induced by multiplication respectively they are isomorphisms. In particular,*

$$\begin{array}{r} \pi_*\mathcal{L}\mathrm{og}^{(k)} \to \pi_*\mathcal{L}\mathrm{og}^{(k-1)} \ , \\ R^1\pi_!\mathcal{L}\mathrm{og}^{(k)} \to R^1\pi_!\mathcal{L}\mathrm{og}^{(k-1)} \end{array}$$

*are zero.*

*Proof.* Both cases are treated at the same time. First consider $R\pi_*$ by induction on $k$: for $k=1$, the assertions follows directly from the definition of $\mathcal{L}og^{(1)}_{G/S}$. For $k \geq 1$, consider the diagram

$$\begin{array}{ccccccccc} 0 \to & \pi^*\mathcal{H}\otimes\pi^*\operatorname{Sym}^k\mathcal{H} & \to & \mathcal{L}og^{(1)}_{G/S}\otimes\pi^*\operatorname{Sym}^k\mathcal{H} & \to & \pi^*\operatorname{Sym}^k\mathcal{H} & \to 0 \\ & \downarrow & & \downarrow & & \downarrow & \\ 0 \to & \pi^*\operatorname{Sym}^{k+1}\mathcal{H} & \to & \mathcal{L}og^{(k+1)}_{G/S} & \to & \mathcal{L}og^{(k)}_{G/S} & \to 0, \end{array}$$

where the first two vertical maps are induced by the multiplication map and the last one is the canonical inclusion. Using lemma A.1.4, the boundary map for $R\pi_*$ induces

$$\begin{array}{ccc} \operatorname{Sym}^k\mathcal{H} & \longrightarrow & \mathcal{H}^\vee\otimes\mathcal{H}\otimes\operatorname{Sym}^k\mathcal{H} \\ \big\downarrow\cong & & \big\downarrow\mathrm{id}\otimes\mathrm{mult} \\ \pi_*\mathcal{L}og^{(k)}_{G/S} & \xrightarrow{\gamma} & \mathcal{H}^\vee\otimes\operatorname{Sym}^{k+1}\mathcal{H} \end{array}$$

The left vertical arrow is an isomorphism by induction. This shows that $\gamma$ is induced by the multiplication map and hence is injective. It follows that $\operatorname{Sym}^{k+1}\mathcal{H}\cong\pi_*\mathcal{L}og^{(k+1)}_{G/S}$. If $G/S$ is a torus, $\gamma$ is even an isomorphism and $R^2\pi_*=0$.

Now consider $R\pi_!$ by induction on $k$: For $k=1$ the boundary map

$$R^1\pi_!\mathbb{Q}_l \to \mathcal{H}(-1)$$

is given as the composition

$$\begin{array}{ccccc} \mathbb{Q}_l\otimes R^1\pi_!\mathbb{Q}_l & \longrightarrow & R^1\pi_*\pi^*\mathcal{H}\otimes R^1\pi_!\mathbb{Q}_l & \xrightarrow{\cup} & R^2\pi_!\pi^*\mathcal{H} \\ \big\downarrow\cong & & \big\downarrow\cong & & \big\downarrow\cong \\ \mathbb{Q}_l\otimes\mathcal{H}(-1) & \xrightarrow{\mathrm{can}\otimes\mathrm{id}} & \mathcal{H}^\vee\otimes\mathcal{H}\otimes\mathcal{H}(-1) & \xrightarrow{\mathrm{mult}\otimes\mathrm{id}} & \mathcal{H}(-1) \end{array}$$

where the first map is induced by the canonical map $\mathbb{Q}_l\to R^1\pi_*\pi^*\mathcal{H}$. This shows that the boundary map is an isomorphism. For $k\geq 1$ the next boundary morphism for $R\pi_!$ in the first diagram of this proof induces

$$\begin{array}{ccc} \mathcal{H}(-1)\otimes\operatorname{Sym}^k\mathcal{H} & \xrightarrow{\sim} & \mathcal{H}(-1)\otimes\operatorname{Sym}^k\mathcal{H} \\ \big\downarrow & & \big\downarrow\mathrm{mult} \\ R^1\pi_!\mathcal{L}og^{(k)}_{G/S} & \xrightarrow{\widetilde{\gamma}} & (\operatorname{Sym}^{k+1}\mathcal{H})(-1). \end{array}$$

Hence $\widetilde{\gamma}$ is surjective, hence an isomorphism for dimension reasons. □

## A.2 The universal property of the logarithm

The main result in this section is the equivalence between unipotent sheaves on $G$ and sheaves on $S$ with an action of a certain symmetric algebra (see theorem A.2.5). This is a formulation of the universal property of the logarithm sheaf.

In order to explain the universal property of the logarithm sheaf, we need the notion of a unipotent sheaf.

**Definition A.2.1.** *A lisse sheaf $\mathcal{F}$ on $G$ is* unipotent of length $n$ *relative to $\pi : G \to S$, if there exists a filtration*

$$\mathcal{F} = A^0\mathcal{F} \supset \ldots \supset A^{n+1}\mathcal{F} = 0$$

*such that there are lisse sheaves $\mathcal{G}_i$ on $S$ with $\mathrm{Gr}^i_A \mathcal{F} \cong \pi^*\mathcal{G}_i$.*

**Definition A.2.2.** *Let $\widehat{\mathfrak{U}}_{\mathcal{H}}$ be the* completion of the universal enveloping algebra *of the abelian Lie algebra $\mathcal{H}$ at the augmentation ideal*

$$I := \ker(\widehat{\mathfrak{U}}_{\mathcal{H}} \to \mathbb{Q}_l).$$

*As $\mathbb{Q}_l$-sheaves on $S$, we have $\widehat{\mathfrak{U}}_{\mathcal{H}} \cong \prod_{k\geq 0} \mathrm{Sym}^k \mathcal{H}$.*

Note that the morphism $\pi^*\mathcal{H} \to \mathcal{L}\mathrm{og}^{(1)}$ induces maps

$$\pi^* \mathrm{Sym}^k \mathcal{H} \otimes \mathcal{L}\mathrm{og}^{(n)} \to \pi^* \mathrm{Sym}^k \mathcal{H} \otimes \mathcal{L}\mathrm{og}^{(n-k)} \to \mathcal{L}\mathrm{og}^{(n)}$$

which give an action of $\widehat{\mathfrak{U}}_{\mathcal{H}}/I^{n+1}$ on $\mathcal{L}\mathrm{og}^{(n)}$.

We want to define a Lie algebra action of the (abelian) Lie algebra $\mathcal{H}$ on $e^*\mathcal{F}$ for all unipotent $\mathcal{F}$.

**Lemma A.2.3.** *(compare [BeL], proof of 1.2.6) Let $\mathcal{F}$ be a unipotent sheaf of length $n$. There is a canonical isomorphism for all $k \geq n$*

$$\pi_*\underline{\mathrm{Hom}}(\mathcal{L}\mathrm{og}^{(k)}, \mathcal{F}) \xrightarrow{\sim} e^*\mathcal{F}$$

*induced by evaluation at $\mathbb{Q}_l \to e^*\mathcal{L}\mathrm{og}^{(n)}$. This induces an action of $\widehat{\mathfrak{U}}_{\mathcal{H}}/I^{n+1}$ on $e^*\mathcal{F}$*

*Proof.* By induction on $n$ we show that the canonical splitting of $e^*\mathcal{L}\mathrm{og}^{(1)} \to e^*\mathbb{Q}_l$ gives an isomorphism

$$e^*\mathcal{F} \cong \pi_*\underline{\mathrm{Hom}}(\mathcal{L}\mathrm{og}^{(k)}, \mathcal{F}).$$

for all $k \geq n$. For $n = 0$ we have for lisse $\mathcal{G}$ on $S$:

$$\pi_* \underline{\mathrm{Hom}}(\mathcal{L}\mathrm{og}^{(k)}, \pi^* \mathcal{G}) \cong \pi_* \underline{\mathrm{Hom}}(\mathcal{L}\mathrm{og}^{(k)}, \pi^! \mathcal{G}(-1)[-2]) \cong \underline{\mathrm{Hom}}(R^2 \pi_! \mathcal{L}\mathrm{og}^{(k)}, \mathcal{G}(-1)) \cong \mathcal{G}$$

according to lemma A.1.4. Assume that the claim is proved for unipotent sheaves of length $n-1$. Then for $\mathcal{F}$ unipotent of length $n$ we have an exact sequence:

$$0 \to A^n \mathcal{F} \to \mathcal{F} \to \mathcal{F}/A^n \mathcal{F} \to 0$$

and $\mathcal{F}/A^n\mathcal{F}$ is unipotent of length $n-1$. Then for $k \geq n-1$

$$\underline{\mathrm{Hom}}(\mathcal{L}\mathrm{og}^{(k)}, \mathcal{F}/A^n \mathcal{F}) \cong e^* \mathcal{F}/A^n \mathcal{F}$$

by induction. The same holds for $A^n \mathcal{F}$. Consider the commutative diagram:

$$\begin{array}{ccccccccc} 0 & \to & \pi_* \underline{\mathrm{Hom}}(\mathcal{L}\mathrm{og}^{(k)}, A^n \mathcal{F}) & \to & \pi_* \underline{\mathrm{Hom}}(\mathcal{L}\mathrm{og}^{(k)}, \mathcal{F}) & \to & \pi_* \underline{\mathrm{Hom}}(\mathcal{L}\mathrm{og}^{(k)}, \mathcal{F}/A^n \mathcal{F}) & \stackrel{\alpha}{\to} & R^1 \pi_* \underline{\mathrm{Hom}}(\mathcal{L}\mathrm{og}^{(k)}, A^n \mathcal{F}) \\ & & \downarrow\cong & & \downarrow & & \downarrow\cong & & \downarrow \\ 0 & \to & e^* A^n \mathcal{F} & \to & e^* \mathcal{F} & \to & e^* \mathcal{F}/A^n \mathcal{F} & \to & 0 \end{array}$$

It suffices to show that the boundary morphism $\alpha : \pi_* \underline{\mathrm{Hom}}(\mathcal{L}\mathrm{og}^{(k)}, \mathcal{F}/A^n \mathcal{F}) \to R^1 \pi_* \underline{\mathrm{Hom}}(\mathcal{L}\mathrm{og}^{(k)}, A^n \mathcal{F})$ is zero. For this consider

$$\begin{array}{ccc} \pi_* \underline{\mathrm{Hom}}(\mathcal{L}\mathrm{og}^{(k-1)}, \mathcal{F}/A^n \mathcal{F}) & \longrightarrow & R^1 \pi_* \underline{\mathrm{Hom}}(\mathcal{L}\mathrm{og}^{(k-1)}, A^n \mathcal{F}) \\ \Big\downarrow \cong & & \Big\downarrow \beta \\ \pi_* \underline{\mathrm{Hom}}(\mathcal{L}\mathrm{og}^{(k)}, \mathcal{F}/A^n \mathcal{F}) & \stackrel{\alpha}{\longrightarrow} & R^1 \pi_* \underline{\mathrm{Hom}}(\mathcal{L}\mathrm{og}^{(k)}, A^n \mathcal{F}) \end{array}$$

and it suffices to prove that $\beta$ is zero. But $A^n \mathcal{F} = \pi^* \mathcal{G} = \pi^! \mathcal{G}(-1)[-2]$ for some $\mathcal{G}$ on $S$, so that

$$R^1 \pi_* R\underline{\mathrm{Hom}}(\mathcal{L}\mathrm{og}^{(k)}, \pi^! \mathcal{G}(-1)[-2]) \cong \underline{\mathrm{Hom}}(R^1 \pi_! \mathcal{L}\mathrm{og}^{(k)}, \mathcal{G}(-1)) \ .$$

It follows that $\beta$ is zero, as the transition maps

$$R^1 \pi_! \mathcal{L}\mathrm{og}^{(k)} \to R^1 \pi_! \mathcal{L}\mathrm{og}^{(k-1)}$$

are zero according to lemma A.1.4 c). The natural pairing for $\underline{\mathrm{Hom}}$ induces an action of $\widehat{\mathfrak{U}}_{\mathcal{H}}$ on $e^* \mathcal{F}$ via

$$e^* \mathcal{L}\mathrm{og} \otimes e^* \mathcal{F} \cong \pi_* \underline{\mathrm{Hom}}(\mathcal{L}\mathrm{og}, \mathcal{L}\mathrm{og}) \otimes \pi_* \underline{\mathrm{Hom}}(\mathcal{L}\mathrm{og}, \mathcal{F}) \to \pi_* \underline{\mathrm{Hom}}(\mathcal{L}\mathrm{og}, \mathcal{F}) \cong e^* \mathcal{F} \ .$$

□

**Corollary A.2.4.** *The action of* $\widehat{\mathfrak{U}}_{\mathcal{H}}/I^{n+1}$ *on*

$$e^*\mathcal{L}\mathrm{og}^{(n)}_{G/S} \cong \widehat{\mathfrak{U}}_{\mathcal{H}}/I^{n+1}$$

*is given by the algebra multiplication in* $\widehat{\mathfrak{U}}_{\mathcal{H}}$.

**Theorem A.2.5.** *(compare: [BeL] 1.2.10 v)) The functor* $\mathcal{F} \longmapsto e^*\mathcal{F}$ *induces an equivalence of the category of unipotent sheaves of length* n *relative* $\pi : G \to S$ *and the category of* $\widehat{\mathfrak{U}}_{\mathcal{H}}/I^{n+1}$*-modules on* $S$.

*Proof.* The functor

$$\mathcal{F}_0 \longmapsto \pi^*\mathcal{F}_0 \otimes_{\widehat{\mathfrak{U}}_{\mathcal{H}}/I^{n+1}} \mathcal{L}\mathrm{og}^{(n)}_{G/S}$$

provides a quasi-inverse. The isomorphism

$$\pi^*e^*\mathcal{F} \otimes_{\widehat{\mathfrak{U}}_{\mathcal{H}}/I^{n+1}} \mathcal{L}\mathrm{og}^{(n)}_{G/S} \cong \mathcal{F}$$

is given by the evaluation map using the description in lemma A.2.3 of $e^*\mathcal{F}$. That it is indeed an isomorphism can be checked after pull-back via $e$. □

**Corollary A.2.6.** *Let* $f : G \to G'$ *be an isogeny, then* $\mathcal{H} \to \mathcal{H}'$ *is an isomorphism and* $\mathcal{L}\mathrm{og}^{(n)}_{G/S} \cong f^*\mathcal{L}\mathrm{og}^{(n)}_{G'/S}$. *In particular,* $\mathcal{L}\mathrm{og}^{(n)}_{G/S}$ *is invariant under translation by torsion sections* $t : S \to G$ *and there is a canonical isomorphism*

$$\varrho_t : t^*\mathcal{L}\mathrm{og}^{(n)}_{G/S} \cong e^*\mathcal{L}\mathrm{og}^{(n)}_{G/S}.$$

*Proof.* The sheaf $\mathcal{L}\mathrm{og}^{(n)}_{G/S}$ is characterized by the action of $\widehat{\mathfrak{U}}_{\mathcal{H}}/I^{n+1}$ on $e^*\mathcal{L}\mathrm{og}^{(n)}_{G/S}$. As $e^*f^* = e'^*$ the claim follows. If $f$ is an isogeny which maps $t$ to $e$, then $\varrho_t$ is given by the composition $t^*\mathcal{L}\mathrm{og}^{(n)}_{G/S} \cong e^*\mathcal{L}\mathrm{og}^{(n)}_{G'/S} \cong e^*\mathcal{L}\mathrm{og}^{(n)}_{G/S}$ and is independent of $f$. □

**Remark:** The above action of $\mathcal{H}$ as a Lie algebra is related to the monodromy action along each fibre on $e^*\mathcal{F}$ as follows: The monodromy gives an action of $\mathcal{H}_s$ on $e^*\mathcal{F}_s$, which is the composition of

$$\mathcal{H} \to \widehat{\mathfrak{U}}_{\mathcal{H}}$$
$$t \mapsto \exp t = \sum_{k \geq 0} \frac{t^i}{i!}$$

and the above action of $\widehat{\mathfrak{U}}_{\mathcal{H}}$ on $e^*\mathcal{F}$.

## A.3 The polylogarithm

Let $U = G \setminus e(S)$ be the complement of the unit section of $G$.

**Lemma A.3.1 ([BeL] 1.3.3).** *Let $\mathcal{L}\mathrm{og}_U$ be the restriction of $\mathcal{L}\mathrm{og}_{G/S}$ to $U$ and $\pi_U : U \to S$ the restriction of $\pi$. Then*

$$\mathrm{Hom}_S(\mathcal{H}, R^1\pi_{U*}\,\mathcal{L}\mathrm{og}_U(1)) = \mathrm{Hom}_S(\mathcal{H}, \mathcal{H}) \ .$$

*Proof.* By purity there is an exact sequence

$$\begin{array}{ccccc} 0 \to & R^1\pi_*\,\mathcal{L}\mathrm{og}(1) & \to & R^1\pi_{U*}\,\mathcal{L}\mathrm{og}_U(1) & \to & e^*\,\mathcal{L}\mathrm{og} \\ & & & & & \| \\ & & & & & \prod_{k\geq 0} \mathrm{Sym}^k\,\mathcal{H} \end{array}$$

and $R^1\pi_*\,\mathcal{L}\mathrm{og}(1)$ is isomorphic to $\mathbb{Q}_l$ if $G/S$ is a torus and equal to 0 if $G/S$ is an elliptic curve (see lemma A.1.4). The lemma follows from this by weight considerations. □

This lemma together with lemma A.1.4 implies that the Leray spectral sequence for $R\pi_{U*}$ induces an isomorphism:

$$\mathrm{Ext}^1_U(\pi_U^*\mathcal{H}, \mathcal{L}\mathrm{og}_U(1)) \cong \mathrm{Hom}_U(\mathcal{H}, \mathcal{H}) \ .$$

**Definition A.3.2.** *Let*

$$\mathcal{P}\mathrm{ol}_{G/S} \in \mathrm{Ext}^1_U((\pi^*\mathcal{H})_U, \mathcal{L}\mathrm{og}_U(1))$$

*be the preimage of the identity in* $\mathrm{Hom}_S(\mathcal{H}, \mathcal{H})$*. This class is called the* (small) polylogarithm extension.

*By abuse of notation we also denote* $\mathcal{P}\mathrm{ol}$ *the (pro)-sheaf representing the extension class (it is unique up to unique isomorphism).*

*For an elliptic curve $E/S$, we call $\mathcal{P}\mathrm{ol}_{E/S}$ the* elliptic polylogarithm.

*For $\mathbb{G}_m/S$, we call $\mathcal{P}\mathrm{ol}_{\mathbb{G}_m/S}$ the* classical or cyclotomic polylogarithm.

By definition we have an exact sequence on $U$

$$0 \to \mathcal{L}\mathrm{og}_U(1) \to \mathcal{P}\mathrm{ol}_{G/S} \to \pi_U^*\mathcal{H} \to 0 \ .$$

It is clear that $\mathcal{P}\mathrm{ol}$ is compatible with base change.

To compare the polylogarithm to the Eisenstein classes (see C.2.2) we need also another version of the above extension.

**Definition A.3.3.** *The* cohomological polylogarithm $\mathcal{P}\mathrm{ol}^{\mathrm{coh}}$ *is the class in*

$$\mathcal{P}\mathrm{ol}^{\mathrm{coh}} \in \mathrm{Ext}^1_U(\mathbb{Q}_l, \pi^*\mathcal{H}^\vee \otimes \mathcal{L}\mathrm{og}_U(1))$$

*obtained by tensoring* $\mathcal{P}\mathrm{ol}$ *with* $\mathcal{H}^\vee$ *and pull-back with the standard map* $\mathbb{Q}_l \to \mathcal{H}^\vee \otimes \mathcal{H}$. *The corresponding sheaf will be denoted by* $\mathcal{P}\mathrm{ol}^{\mathrm{coh}}$.

**Remark:** We call this class cohomological because it appears as an element of cohomology of appropriate varieties.

# B Degeneration in the local situation

In this appendix we study the degeneration of the elliptic polylog into the classical polylog induced by the degeneration of an elliptic curve into a torus. The main results are given in theorem B.1.3. It is enough to consider the situation locally around a point of bad reduction.

## B.1 Degeneration of $\mathcal{L}\mathrm{og}$ and $\mathcal{P}\mathrm{ol}$

Let $S$ be the spectrum of a complete discrete valuation ring with quotient field $K$ and residue field $\kappa$. We put $\eta := \mathrm{Spec}\, K$ and $s := \mathrm{Spec}\, \kappa$. Let $E/\eta$ be an elliptic curve with level-$N$-structure. Let $\widetilde{E}$ its Néron model. Denote by $e : S \to \widetilde{E}$ the identity section. We consider the cartesian diagram

$$\begin{array}{ccccc} E & \xrightarrow{j} & \widetilde{E} & \xleftarrow{i} & \widetilde{E}_s \\ {\scriptstyle\pi}\downarrow & & {\scriptstyle\widetilde{\pi}}\downarrow & & \downarrow{\scriptstyle\widetilde{\pi}} \\ \eta & \xrightarrow{j'} & S & \xleftarrow{i'} & s \end{array}$$

Let $\widetilde{E}^0$ be the connected component of the identity section and *assume* that there is an isomorphism

$$\widetilde{E}_s = \mathbb{G}_{m,s} \times \mathbb{Z}/N\mathbb{Z}$$

with $\widetilde{E}^0 = \mathbb{G}_{m,s} \times \{0\}$. We define

$$\mathbb{G}_m^{(v)} := \mathbb{G}_{m,s} \times \{v\}$$

and let

$$e_v : S \to \widetilde{E}$$

be a global section which specialize to the identity section of $\mathbb{G}_m^{(v)}$. We are going to study the behaviour of the logarithmic and the polylogarithmic sheaves under the functor $i^* j_*$.

**Remark:** We can for example put $S$ the completion of $\overline{M}$ at $\infty$. Then $E$ and $\widetilde{E}$ are the base changes of what was called $E$ and $\widetilde{E}$ in section 1.1.

By $\mathcal{H}_\eta$ and $\mathcal{H}_s$ we denote the sheaves $\mathcal{H}$ (see definition A.1.1) for $E/\eta$ respectively $\widetilde{E}^0_s/s$.

**Lemma B.1.1.** *There is a canonical isomorphism $\mathbb{Q}_l(1) = \mathcal{H}_s \to i'^* j'_* \mathcal{H}_\eta$. Let*

$$\mathbb{Q}_l(1) \xrightarrow{\partial} \mathcal{H}_\eta$$

*be its adjoint. Let*

$$p(1) : \mathbb{Q}_l(1) \otimes \mathcal{H}_\eta \xrightarrow{\partial \otimes \mathrm{id}} \mathcal{H}_\eta \otimes \mathcal{H}_\eta \to \mathbb{Q}_l(1)$$

*be given by the intersection pairing. We get a short exact sequence*

$$0 \to \mathbb{Q}_l(1) \xrightarrow{\partial} \mathcal{H}_\eta \xrightarrow{p} \mathbb{Q}_l \to 0 \ .$$

*Proof.* $\mathcal{H}_\eta$ has non-trivial monodromy. Hence the monodromy filtration $M_*$ has two non-trivial steps and we have the monodromy isomorphism(([Del1] 1.6.14)

$$\mathrm{Gr}_1 \mathcal{H}_\eta(1) \to \mathrm{Gr}_{-1} \mathcal{H}_\eta \ .$$

$M_{-1}\mathcal{H}_\eta$ extends to a lisse sheaf on $S$. It is enough to show that $i'^* j'_* M_{-1}\mathcal{H}_\eta = \mathcal{H}_s$. We sketch a proof of this well-known fact: Let $\overline{E}$ be the compactification of $\widetilde{E}$ by a generalized elliptic curve. $\overline{E}_s$ is the proper and singular Néron-$N$-gon. Consider the commutative diagram

$$\begin{array}{ccc} E & \xrightarrow{\bar{j}} & \overline{E} \\ \pi \downarrow & & \downarrow \overline{\pi} \\ \eta & \xrightarrow{j'} & S \end{array}$$

Note that $R^1\overline{\pi}_*\mathbb{Q}_l(1) = R^1\widetilde{\pi}_!\mathbb{Q}_l(1) = (R^1\widetilde{\pi}_*\mathbb{Q}_l)^\vee$. It is isomorphic to $\mathcal{H}_\eta$ and $\mathcal{H}_s$ over $\eta$ and $s$ respectively. We consider the Leray spectral sequences for the compositions $\overline{\pi} \circ \overline{j}$ and $j' \circ \pi$ which converge to the same thing. The five term sequences give

$$0 \to i'^* R^1 \overline{\pi}_* \mathbb{Q}_l(1) \to i'^* R^1 (\overline{\pi} \circ \overline{j})_* \mathbb{Q}_l(1) \to i'^* \overline{\pi}_* \mathbb{Q}_l \to \cdots$$

respectively

$$0 \to i'^* \mathbb{Q}_l \to i'^* R^1 (j' \circ \pi_\eta)_* \mathbb{Q}_l(1) \to i'^* j'_* R^1 \pi_{\eta *} \mathbb{Q}_l(1) \to 0$$

But $i'^* R^1 \overline{\pi}_* \mathbb{Q}_l(1) = \mathbb{Q}_l(1)$ and comparison of the two sequences gives the result. □

The consequences of this lemma for $i'^* R^n j'_* \operatorname{Sym}^n \mathcal{H}_\eta$ have already been drawn in 2.1.3.

For reference we also note:

**Corollary B.1.2.** *There are canonical isomorphisms:*

a) $$i'^* R^1 j'_* e^* \mathcal{L}\mathrm{og}_{E/\eta}(1) \cong \prod_{k \geq 0} \mathbb{Q}_l \ ,$$

b) $$i'^* R^1 j'_* \left( \mathcal{H}^\vee \otimes e^* \mathcal{L}\mathrm{og}_{E/\eta}(1) \right) \cong \prod_{k \geq 1} \mathbb{Q}_l \oplus \prod_{k \geq 0} \mathbb{Q}_l(-1) \ .$$

*Proof.* We have $e^* \mathcal{L}\mathrm{og}_{E/\eta}(1) \cong \prod_{k \geq 0} \operatorname{Sym}^k \mathcal{H}_\eta(1)$ and a canonical isomorphism

$$\mathcal{H}_\eta^\vee \otimes \prod_{k \geq 0} \operatorname{Sym}^k \mathcal{H}_\eta(1) \cong \prod_{k \geq 0} \left( \operatorname{Sym}^{k-1} \mathcal{H}_\eta(1) \oplus \operatorname{Sym}^{k+1} \mathcal{H}_\eta \right).$$

□

**Theorem B.1.3 (Degeneration, cf. [BeL] 1.5). a)** *There is a canonical isomorphism*

$$i^* j_* \mathcal{L}\mathrm{og}_{E/\eta} |_{\mathbb{G}_m^{(v)}} \cong \mathcal{L}\mathrm{og}_{\mathbb{G}_m^{(v)}/s}$$

**b)** $\mathcal{P}\mathrm{ol}_{E/\eta} \in \operatorname{Ext}^1_{E \smallsetminus e(\eta)}(\pi^* \mathcal{H}_\eta, \mathcal{L}\mathrm{og}_{E/\eta}(1))$ *has residue zero, i.e., the sequence*

$$0 \to i^* j_* \mathcal{L}\mathrm{og}_{E/\eta}(1) \to i^* j_* \mathcal{P}\mathrm{ol}_{E/\eta} \to \mathbb{Q}_l(1) \to 0$$

*is still exact. On $\widetilde{E}^0_s \smallsetminus e(s) = \mathbb{G}_m^{(0)} \smallsetminus e(s)$ there is a canonical isomorphism*

$$i^* j_* \mathcal{P}\mathrm{ol}_{E/\eta}|_{\mathbb{G}_m^{(0)}} \cong \mathcal{P}\mathrm{ol}_{\mathbb{G}_m} .$$

*On $\mathbb{G}_m^{(v)}$ for $v \neq 0$ the extension $i^* j_* \mathcal{P}\mathrm{ol}_{E/\eta}|_{\mathbb{G}_m^{(v)}}$ splits, more precisely*

$$\mathrm{Ext}^1_{\mathbb{G}_m^{(v)}}(\mathbb{Q}_l(1), i^* j_* \mathcal{L}\mathrm{og}_{E/\eta}(1)) = 0.$$

The theorem will be proved in the next section.

**Corollary B.1.4.** *Let $t \in \widetilde{E}(S) \smallsetminus e(S)$. Then $t^* \mathcal{P}\mathrm{ol}_{E/\eta} \in \mathrm{Ext}^1_\eta(\mathcal{H}_\eta, t^* \mathcal{L}\mathrm{og}_{E/\eta}(1))$ has residue zero, i.e.,*

$$0 \to i'^* j'_* t^* \mathcal{L}\mathrm{og}_{E/\eta}(1) \to i'^* j'_* t^* \mathcal{P}\mathrm{ol}_{E/\eta} \to \mathbb{Q}_l(1) \to 0$$

*is still exact. For $t \in \widetilde{E}^0(S) \smallsetminus e(S)$, there is a canonical isomorphism*

$$i'^* j'_* t^* \mathcal{P}\mathrm{ol}_{E/\eta} \cong t^* \mathcal{P}\mathrm{ol}_{\mathbb{G}_m^{(0)}}$$

*and for $t \in \widetilde{E}(S) \smallsetminus \widetilde{E}^0(S)$ the extension*

$$0 \to i'^* j'_* t^* \mathcal{L}\mathrm{og}_{E/\eta}(1) \to i'^* j'_* t^* \mathcal{P}\mathrm{ol}_{E/\eta} \to i'^* j'_* \mathcal{H}_\eta \to 0$$

*splits.*

*Proof.* This follows immediately from part b) of the theorem using lemma B.2.1 below. □

**Remark:** One point of the theorem is that the the extensions $\mathcal{P}\mathrm{ol}_{E/\eta}$ and hence $t^* \mathcal{P}\mathrm{ol}_{E/\eta}$ have residue zero under $i^* j_*$. It is a key point that this is *false* for

$$0 \to \mathcal{L}\mathrm{og}_{E/\eta} \otimes \mathcal{H}_\eta^\vee \to \mathcal{P}\mathrm{ol}^{\mathrm{coh}}_{E/\eta} \to \mathbb{Q}_l \to 0 .$$

Indeed, the residue of this sequence will be computed at least after pull-back along torsion sections in theorem C.1.1. Taking certain linear combinations of pull-backs along torsion section this residue can again be zero and is then connected to the residue of the elliptic polylog (see 2.3.2) and hence to the classical polylog.

## B.2 Proof of the degeneration theorem B.1.3

The following lemma allows us to compute the degeneration on the base:

**Lemma B.2.1.** *Let $t : S \to \widetilde{E}$ be a section, then for all $q \geq 0$ and unipotent $l$-adic sheaves $\mathcal{F}$ on $E$*

$$i'^* R^q j'_* t^* \mathcal{F} \cong t^* i^* R^q j_* \mathcal{F}$$

*The same holds for open subschemes $U \subset \widetilde{E}$ and sections $t : S \to U$.*

*Proof.* Adjunction for $j$ induces a map from the right hand side to the left. By induction on the length of the filtration on $\mathcal{F}$, it is enough to show the assertion for $\mathcal{F} = \pi^* \mathcal{G}$. In this case

$$Rj'_* t^* \pi^* \mathcal{G} = t^* \widetilde{\pi}^* Rj'_* \mathcal{G} = t^* Rj_* \pi^* \mathcal{G}$$

where the first equality follows from $\pi \circ t = \mathrm{id}$ and the second is smooth base change. □

**Lemma B.2.2.** *Let $A$ be the canonical filtration on $\mathcal{L}\mathrm{og}$.*

**a)** *The sheaf $i^* j_* \mathcal{L}\mathrm{og}^{(n)}_{E/\eta}$ is unipotent with $\mathrm{Gr}^k_{i^* j_* A}\, i^* j_* \mathcal{L}\mathrm{og}^{(n)}_{E/\eta} \cong \mathbb{Q}_l(k)$ for $0 \leq k \leq n$.*

**b)** *The sheaves $i^* R^1 j_* \mathcal{L}\mathrm{og}^{(n)}_{E/\eta}$ are unipotent with ($0 \leq k \leq n$)*

$$\mathrm{Gr}^k_{i^* R^1 j_* A}\, i^* R^1 j_* \mathcal{L}\mathrm{og}^{(n)}_{E/\eta} \cong \mathbb{Q}_l(-1)\ .$$

*Proof.* It is enough to consider the case when $\kappa$ is finitely generated over its prime field. The general case follows by base change. Hence we can assume that there are no non-trivial morphisms between $\mathbb{Q}_l(k)'s$ for different $k$.

We prove a) and b) together by induction on $n$: for $n = 0$, $\mathcal{L}\mathrm{og}^{(0)} = \mathbb{Q}_l$ and there is nothing to prove. For $n > 0$ consider

$$0 \to \pi^* \operatorname{Sym}^n \mathcal{H}_\eta \to \mathcal{L}\mathrm{og}^{(n)}_{E/\eta} \to \mathcal{L}\mathrm{og}^{(n-1)}_{E/\eta} \to 0.$$

By induction and 2.1.3 the boundary map

$$i^* j_* \mathcal{L}\mathrm{og}^{(n-1)}_{E/\eta} \to i^* R^1 j_* \operatorname{Sym}^n \mathcal{H}_\eta \cong \mathbb{Q}_l(-1)$$

vanishes for weight reasons. As $i^* j_* \operatorname{Sym}^n \mathcal{H}_\eta \cong \mathbb{Q}_l(n)$ the claim follows. □

Let $e_v : S \to \widetilde{E}$ be the section defined in B.1. By the characterization of unipotent sheaves in theorem A.2.5, it remains to prove that $e_v^* i^* j_* \mathcal{L}\mathrm{og}^{(n)}$ is isomorphic to $\widehat{\mathfrak{U}}_{\mathcal{H}_s}/I^{n+1}$ with its canonical $\widehat{\mathfrak{U}}_{\mathcal{H}_s}$-module structure.

**Lemma B.2.3.** *There is a canonical isomorphism of $\widehat{\mathfrak{U}}_{\mathcal{H}_s}$-modules*

$$e_v^* i^* j_* \,\mathcal{L}\mathrm{og} \cong \widehat{\mathfrak{U}}_{\mathcal{H}_s}$$

*Proof.* We have an isomorphism of $\widehat{\mathfrak{U}}_{\mathcal{H}_s}$-modules

$$e_v^* i^* j_* \,\mathcal{L}\mathrm{og} \cong i'^* j'_* e_v^* \,\mathcal{L}\mathrm{og} \cong i'^* j'_* \widehat{\mathfrak{U}}_{\mathcal{H}} \cong \widehat{\mathfrak{U}}_{\mathcal{H}_s}$$

by B.2.1 for the first and 2.1.3 for the last isomorphism. □

*Proof.* (of theorem B.1.3). The lemma proves part a) of the theorem.

To prove part b) for $v \neq 0$ we remark that $i^* j_* \,\mathcal{L}\mathrm{og}_{E/\eta}|_{\mathbb{G}_m^{(v)}} \cong \mathcal{L}\mathrm{og}_{\mathbb{G}_m^{(v)}}$ by a) and that

$$\pi_* \,\mathcal{L}\mathrm{og}_{\mathbb{G}_m^{(v)}} = 0 \text{ and } R^1\pi_* \,\mathcal{L}\mathrm{og}_{\mathbb{G}_m^{(v)}} \cong \mathbb{Q}_l(-1)$$

by lemma A.1.4. The claim follows from the Leray spectral sequence.

To prove part b) for $v = 0$, consider res : $\mathrm{Ext}^1_{U_\eta}(\mathcal{L}\mathrm{og}_{E/\eta}(1), \pi^*\mathcal{H}_\eta) \to \mathrm{Hom}_{U_s^0}(i^* j_* \mathcal{H}_\eta, i^* R^1 j_* \,\mathcal{L}\mathrm{og}_{E/\eta}(1))$ with $U_s^0 = \widetilde{E}_s^0 \smallsetminus e(s)$. Note that $\mathrm{res}(\mathcal{P}\mathrm{ol}_{E/\eta})$ vanishes by corollary B.2.2 for weight reasons. (As in the proof of B.2.2 we can assume without loss of generality that $\kappa$ is finitely generated over its prime field.) Hence we have a short exact sequence

$$0 \to i^* j_* \,\mathcal{L}\mathrm{og}(1)_{E/\eta} \to i^* j_* \,\mathcal{P}\mathrm{ol}_{E/\eta} \to i^* j_* \pi^* \mathcal{H}_\eta \to 0\ .$$

To identify this extension on $\mathbb{G}_m^{(0)}$, consider the commutative diagram

$$\begin{array}{ccc}
\mathrm{Ext}^1_{U_\eta}(\pi^*\mathcal{H}_\eta, \mathcal{L}\mathrm{og}_{E/\eta}(1))^{\mathrm{res}=0} & \longrightarrow & \mathrm{Hom}_\eta(\mathcal{H}_\eta, e^* \,\mathcal{L}\mathrm{og}_{E/\eta}) \\
{\scriptstyle i^* j_*}\downarrow & & \downarrow{\scriptstyle i'^* j'_*} \\
\mathrm{Ext}^1_{U_s^0}(i^* j_* \pi^* \mathcal{H}_\eta, i^* j_* \,\mathcal{L}\mathrm{og}_{E/\eta}(1)) & \longrightarrow & \mathrm{Hom}_s(i'^* j'_* \mathcal{H}_\eta, i'^* j'_* e^* \,\mathcal{L}\mathrm{og}_{E/\eta}) \\
\cong\downarrow & & \downarrow\cong \\
\mathrm{Ext}^1_{U_s^0}(\widetilde{\pi}_s^* \mathcal{H}_s, \mathcal{L}\mathrm{og}_{\mathbb{G}_m^{(0)}}(1)) & \longrightarrow & \mathrm{Hom}_s(\mathcal{H}_s, e^* \,\mathcal{L}\mathrm{og}_{\mathbb{G}_m^{(0)}})
\end{array}$$

Hence the class of $i^* j_* \,\mathcal{P}\mathrm{ol}_{E/\eta}$ is mapped to the identity. This proves part b) because the polylog was uniquely determined by its residue. □

# C The comparison theorem and values at torsion sections

We are going to show that the pull-back of the elliptic polylog at torsion sections is given by the Eisenstein symbol. This is done by showing that the Eisenstein symbol and the pull-back of the polylogarithm at torsion sections both give a section of a certain residue map. As this section is uniquely determined for weight reasons, it is enough to compute the residues of the pull-backs of the polylogarithm and to see that this agrees with the residue of the Eisenstein symbol.

## C.1 Residue of the elliptic polylog at torsion sections

In this section, we work in the universal situation of section 1.1. The $\mathbb{Q}(\mu_N)$-scheme $M$ is the modular curve of elliptic curves with full level-$N$-structure, $E$ is the universal elliptic curve over $M$ and $\widetilde{E}$ is the Néron model over the compactification $\overline{M}$. For further notation we refer to section 1.1.

We are going to use the residue (in the sense of 0.0.1)

$$\begin{aligned}\operatorname{Ext}^1_M(\mathbb{Q}_l, \mathcal{H}^\vee \otimes e^* \mathcal{L}\mathrm{og}_{E/M}(1)) &\to \operatorname{Hom}_{\mathrm{Cusp}}(\mathbb{Q}_l, i'^* R^1 j'_* (\mathcal{H}^\vee \otimes e^* \mathcal{L}\mathrm{og}_{E/M}(1))) \\ &\xrightarrow{2.2.3} \operatorname{Hom}_{\mathrm{Cusp}}(\mathbb{Q}_l, i'^* R^1 j'_* (e^* \mathcal{L}\mathrm{og}_{E/M}(1)))\end{aligned}$$

We compose this with the pull-back to Isom (see section 1.1)

$$\begin{aligned}&\operatorname{Hom}_{\mathrm{Cusp}}(\mathbb{Q}_l, i'^* R^1 j'_* (e^* \mathcal{L}\mathrm{og}_{E/M}(1))) \longrightarrow \\ &\qquad\qquad \operatorname{Hom}_{\mathrm{Isom}}(\mathbb{Q}_l, i'^* R^1 j'_* (e^* \mathcal{L}\mathrm{og}_{E/M}(1))) \ .\end{aligned}$$

By corollary B.1.2 (which applies because we have an isomorphism $\widetilde{E}_{\mathrm{Isom}} \cong \mathbb{G}_m \times \mathbb{Z}/N \times \mathrm{Isom}$) the target space is isomorphic to

$$\operatorname{Hom}_{\mathrm{Isom}}(\mathbb{Q}_l, i'^* R^1 j'_* (e^* \mathcal{L}\mathrm{og}_{E/M}(1))) \cong \operatorname{Hom}_{\mathrm{Isom}}(\mathbb{Q}_l, \prod_{k \geq 1} \mathbb{Q}_l) \ .$$

If we parameterize Isom as in section 1.1, we can identify $\operatorname{Hom}_{\mathrm{Isom}}(\mathbb{Q}_l, \prod_{k\geq 1} \mathbb{Q}_l)$ with

$$\prod_{g \in P(\mathbb{Z}/N) \backslash \mathrm{Gl}_2(\mathbb{Z}/N)} \operatorname{Hom}_{\overline{\mathbb{Q}}}(\mathbb{Q}_l, \prod_{k \geq 1} \mathbb{Q}_l).$$

All in all we consider res as a map

$$\mathrm{res} : \operatorname{Ext}^1_M(\mathbb{Q}_l, \mathcal{H}^\vee \otimes e^* \mathcal{L}\mathrm{og}_{E/M}(1)) \longrightarrow \prod_{g \in P(\mathbb{Z}/N) \backslash \mathrm{Gl}_2(\mathbb{Z}/N)} \prod_{k \geq 1} \mathbb{Q}_l \ .$$

Finally recall that for a torsion section $t \in E[N](M)$ with $t \neq 0$ we have the isomorphism A.2.6

$$\varrho_t : \mathrm{Ext}^1_M(\mathbb{Q}_l, \mathcal{H}^\vee \otimes t^* \mathcal{L}\mathrm{og}_{E/M}(1)) \to \mathrm{Ext}^1_M(\mathbb{Q}_l, \mathcal{H}^\vee \otimes e^* \mathcal{L}\mathrm{og}_{E/M}(1)).$$

We now describe the residue of $\varrho_t(t^* \mathcal{P}\mathrm{ol}^{\mathrm{coh}}_{E/M})$. Recall that $B_k(x)$ is the Bernoulli polynomial as in section 1.2. We identify $x \in \mathbb{R}/\mathbb{Z}$ with its representative in $[0, 1)$.

**Theorem C.1.1 (Residue of values at torsion sections).** *Let $k \geq 0$, $N \geq 3$ and let $t : M \to \widetilde{E}$ be a non zero $N$-torsion section. Then for $g \in P(\mathbb{Z}/N) \setminus \mathrm{Gl}_2(\mathbb{Z}/N) = \mathrm{Isom}(\overline{\mathbb{Q}})$*

$$\mathrm{res}\left(\varrho_t(t^* \mathcal{P}\mathrm{ol}^{\mathrm{coh}}_{E/M})\right)^{k+1}(g) = \frac{-N}{(k+2)k!} B_{k+2}\left(\frac{(gt)_2}{N}\right)$$

*where we have written $gt = ((gt)_1, (gt)_2) \in E[N] \cong (\mathbb{Z}/N)^2$.*

This is due to Beilinson and Levin [BeL] prop. 2.2.3. A different proof can be found in [W] cor. III 3.26. The formula given there differs by a factor $N$ from ours (i.e. the one in [BeL]). This seems to come from a different normalization of the residue. For the convenience of the reader and because we needed to fix the normalizations we give a third proof (also due to Beilinson and Levin from an earlier version of [BeL]) in the last section C.3. The reader is advised to skip it on a first reading.

## C.2 The polylog at torsion sections

We continue in the modular situation of the last section. We want to relate $\varrho_t(t^* \mathcal{P}\mathrm{ol}^{\mathrm{coh}}_{E/M})$ and the Eisenstein symbol.

Recall (cf. section 1.2) that $\psi \in \mathbb{Q}[E[N]]$ is mapped to an element in $\mathbb{Q}[\mathrm{Isom}]$ via the horospherical map $\varrho^k$ (1.2.3) and

$$\mathrm{Eis}^k_l(\varrho^k \psi) \in H^1_{\mathrm{et}}(M, \mathrm{Sym}^k \mathcal{H}(1)) .$$

For convenience, we make the following definition.

**Definition C.2.1.** *For $\psi \in \mathbb{Q}[E[N] \smallsetminus 0]$ put*

$$\psi^* \mathcal{P}\mathrm{ol}^{\mathrm{coh}}_{E/M} = \sum_{t \in E[N]} \psi(t) \varrho_t(t^* \mathcal{P}\mathrm{ol}^{\mathrm{coh}}_{E/M})$$

*as element of* $\mathrm{Ext}^1_M(\mathbb{Q}_l, \mathcal{H}^\vee \otimes e^* \mathcal{L}\mathrm{og}_{E/M}(1))$.

**Theorem C.2.2 (Values at torsion sections).**
*Let $k > 0$, $N \geq 3$ and $\psi \in \mathbb{Q}[E[N] \smallsetminus 0]$, then*

$$\mathrm{Eis}_l^k(\varrho^k \psi) = -N^{k-1}\,\mathrm{pr}\left(\psi^* \,\mathcal{P}\mathrm{ol}^{\mathrm{coh}}\right)^{k+1}$$

*as elements of*

$$\mathrm{Ext}^1_M(\mathbb{Q}_l, \mathrm{Sym}^k \mathcal{H}(1))$$

*with the projection* pr *in lemma 2.2.3.*

**Remark:** The theorem for the universal elliptic curve immediately implies the same theorem for elliptic curves over any base scheme of characteristic zero (which is of finite type over a regular 0- or 1-dimensional scheme) by base-change from the modular situation.

The proof the theorem will take up the rest of this section.

We start with a compatibility result for our different residue maps.

**Lemma C.2.3.** *The following diagram commutes*

$$\begin{array}{ccccc}
H^{k+1}_{\mathrm{et}}(E^k, \mathbb{Q}_l(k+1))(\varepsilon) & \xrightarrow[0.0.1a)]{\mathrm{res}} & H^k_{\mathrm{et}}((\widetilde{E}^0)^k_{\mathrm{Cusp}}, \mathbb{Q}_l(k))(\varepsilon) & \longrightarrow & H^0_{\mathrm{et}}(\mathrm{Isom}, \mathbb{Q}_l(0)) \\
\downarrow \cong & & \downarrow \cong & & \downarrow = \\
H^1_{\mathrm{et}}(M, \mathrm{Sym}^k \mathcal{H}(1)) & \xrightarrow[0.0.1b)]{\mathrm{res}} & H^0_{\mathrm{et}}(\mathrm{Cusp}, i^* R^1 j_* \mathrm{Sym}^k \mathcal{H}) & \longrightarrow & H^0_{\mathrm{et}}(\mathrm{Isom}, \mathbb{Q}_l)\ .
\end{array}$$

*Proof.* Consider the commutative diagram

$$\begin{array}{ccc}
H^{k+1}_{\mathrm{et}}((\widetilde{E}^0)^k, Rj_*\mathbb{Q}_l(k+1))(\varepsilon) & \longrightarrow & H^{k+2}_{\mathrm{et}}((\widetilde{E}^0)^k_{\mathrm{Cusp}}, Ri^!\mathbb{Q}_l(k+1))(\varepsilon) \\
\downarrow \cong & & \downarrow \cong \\
H^{k+1}_{\mathrm{et}}(\overline{M}, R\widetilde{\pi}_* Rj_*\mathbb{Q}_l(k+1))(\varepsilon) & \longrightarrow & H^{k+2}_{\mathrm{et}}(\mathrm{Cusp}, R\widetilde{\pi}_* Ri^!\mathbb{Q}_l(k+1))(\varepsilon) \\
\downarrow \cong & & \downarrow \cong \\
H^{k+1}_{\mathrm{et}}(\overline{M}, Rj'_* R\pi_*\mathbb{Q}_l(k+1))(\varepsilon) & \longrightarrow & H^{k+2}_{\mathrm{et}}(\mathrm{Cusp}, Ri'^! R\widetilde{\pi}_*\mathbb{Q}_l(k+1))(\varepsilon) \\
\downarrow \cong & & \downarrow \cong \\
H^1_{\mathrm{et}}(\overline{M}, Rj'_* \mathrm{Sym}^k \mathcal{H}(1)) & \longrightarrow & H^2_{\mathrm{et}}(\mathrm{Cusp}, Ri'^! \mathrm{Sym}^k(R^1\widetilde{\pi}_*\mathbb{Q}_l(1))(1))
\end{array}$$

Here we have used the dual of proper base change $R\widetilde{\pi}_* i^! \cong i'^! R\widetilde{\pi}_*$ in the middle square. The composition is the residue. For the last line the distinguished triangle

$$i'_* Ri'^! \mathrm{Sym}^k(R^1\widetilde{\pi}_*\mathbb{Q}_l(1)) \to \mathrm{Sym}^k(R^1\widetilde{\pi}_*\mathbb{Q}_l(1)) \to Rj'_* j'^* \mathrm{Sym}^k(R^1\widetilde{\pi}_*\mathbb{Q}_l(1))$$

shows that

$$i'^{*}R^1 j'_* \operatorname{Sym}^k \mathcal{H} \cong i'^{*}R^1 j'_* \operatorname{Sym}^k(R^1\pi_*\mathbb{Q}_l(1)) \cong R^2 i'^{!} \operatorname{Sym}^k(R^1\widetilde{\pi}_*\mathbb{Q}_l(1)) \ .$$

Together with purity in the first line of the diagram, we obtain

$$\begin{array}{ccc} H^{k+1}_{\mathrm{et}}((\widetilde{E}^0)^k, Rj_*\mathbb{Q}_l(k+1))(\varepsilon) & \xrightarrow{\mathrm{res}} & H^{k+2}_{\mathrm{et}}((\widetilde{E}^0)^k_{\mathrm{Cusp}}, \mathbb{Q}_l(k))(\varepsilon) \\ \downarrow\cong & & \downarrow\cong \\ H^1_{\mathrm{et}}(\overline{M}, Rj'_* \operatorname{Sym}^k \mathcal{H}(1)) & \xrightarrow{\mathrm{res}} & H^0_{\mathrm{et}}(\mathrm{Cusp}, i'^{*}R^1 j'_* \operatorname{Sym}^k \mathcal{H}) \end{array}$$

which proves our claim. The second diagram is obtained by first applying base change to Isom and then using the canonical isomorphism $\widetilde{E}^0_{\mathrm{Isom}} = \mathbb{G}_m \times \mathrm{Isom}$. □

*Proof.* (of Theorem C.2.2) We have to compare two elements in $H^1_{\mathrm{et}}(M, \operatorname{Sym}^k \mathcal{H}(1))$ given by the $k+1$-st component of $\psi^* \mathcal{P}\mathrm{ol}^{\mathrm{coh}}$ respectively by $\mathrm{Eis}^k_l(\psi)$.

Let $a: M \to B$ be the structural map. Consider the short exact sequence induced by the Leray spectral sequence

$$0 \to H^1_{\mathrm{et}}(B, a_* \operatorname{Sym}^k \mathcal{H}(1)) \to H^1_{\mathrm{et}}(M, \operatorname{Sym}^k \mathcal{H}(1)) \to H^0_{\mathrm{et}}(B, R^1 a_* \operatorname{Sym}^k \mathcal{H}(1)) \to 0 \ .$$

As $\operatorname{Sym}^k \mathcal{H}$ has no global sections for $k > 0$, the sheaf $a_* \operatorname{Sym}^k \mathcal{H}$ is zero. Hence it is enough to show that our two elements have the same image in $H^1_{\mathrm{et}}(M \times \overline{\mathbb{Q}}, \operatorname{Sym}^k \mathcal{H}(1))$, the vector space underlying $R^1 a_* \operatorname{Sym}^k \mathcal{H}(1)$. Now consider the residue sequence

$$\begin{aligned} 0 \to H^1_{\mathrm{et}}(\overline{M} \times \overline{\mathbb{Q}}, j'_* \operatorname{Sym}^k \mathcal{H}(1)) \to H^1_{\mathrm{et}}(M \times \overline{\mathbb{Q}}, \operatorname{Sym}^k \mathcal{H}(1)) \to \\ H^0_{\mathrm{et}}(\mathrm{Cusp} \times \overline{\mathbb{Q}}, i'^{*}R^1 j'_* \operatorname{Sym}^k \mathcal{H}(1)) \to 0. \end{aligned}$$

The cuspidal cohomology on the left is an intersection cohomology group and hence pure of weight $-k-1$. The term on the right is pure of weight 0. Our elements are Galois-invariant and hence uniquely determined by their image in

$$H^0_{\mathrm{et}}(\mathrm{Cusp} \times \overline{\mathbb{Q}}, i'^{*}R^1 j'_* \operatorname{Sym}^k \mathcal{H}(1)) \subset H^0_{\mathrm{et}}(\mathrm{Isom} \times \overline{\mathbb{Q}}, \mathbb{Q}_l) \ .$$

The commutative diagram obtained in lemma C.2.3 allows us to compare the residue of the $k+1$-st component of $\psi^* \mathcal{P}\mathrm{ol}^{\mathrm{coh}}$ and that of $\mathrm{Eis}^k_l(\psi)$. Hence it is enough to see that the two maps from $\mathbb{Q}[E[N] \smallsetminus 0]$ to $H^0_{\mathrm{et}}(\mathrm{Isom} \times \overline{\mathbb{Q}}, \mathbb{Q}_l)$ given by

$$\begin{aligned} \psi &\mapsto -N^{k-1} \operatorname{res} \operatorname{pr} \left(\psi^* \mathcal{P}\mathrm{ol}^{\mathrm{coh}}\right)^{k+1} \\ \psi &\mapsto \operatorname{res} \circ \mathrm{Eis}_l(\varrho^k \psi) \end{aligned}$$

agree. This follows from theorem C.1.1 and proposition 1.2.1. □

## C.3 The residue computation

In this section we give an exposition of a proof of C.1.1. We follow very closely the arguments of Beilinson and Levin in an earlier preprint version of [BeL].

Because of $\mathrm{Gl}_2$-equivariance it is enough to consider $g = \mathrm{id}$, i.e., residue at the cusp $\infty$. Moreover, we use base change to $\mathbb{C}$ and consider the residue of the underlying local systems. Then the question is local around $\infty$. Hence our base is

$$N\mathbb{Z} \setminus \mathfrak{H} \xrightarrow{\tau \mapsto \exp(2\pi i \frac{\tau}{N})} \mathbb{C}^*$$

The elliptic curve $E_\tau$ is given by $\mathbb{C}/\mathbb{Z} + \mathbb{Z}\tau$ with level-$N$-structure generated by $\frac{1}{N}$, $\frac{\tau}{N}$. Let $t = \frac{a}{N}\tau + \frac{b}{N}$ be a nonzero $N$-torsion point of $E_\tau$. We assume $0 \leq a, b < N$. The fundamental group of the base is $\Gamma_\infty \subset \mathrm{SL}_2(\mathbb{Z})$ generated by $\widetilde{T} = \begin{pmatrix} 1 & N \\ 0 & 1 \end{pmatrix}$ with its natural operation on $\mathfrak{H}$ from the left. We consider $\mathcal{H}$, $t^* \mathcal{L}\mathrm{og}$, $t^* \mathcal{P}\mathrm{ol}$ etc. (we drop the index $E(\mathbb{C})/M(\mathbb{C})$) as local systems of $\mathbb{Q}$-vector spaces on $N\mathbb{Z} \setminus \mathfrak{H}$, i.e., as representations of $\Gamma_\infty$.

**Lemma C.3.1.** *Let $e_1, e_2$ be the classes in $\mathcal{H}_\tau = H_1(E_\tau, \mathbb{Q})$ represented by the loops $[t, t+1], [t, t+\tau]$ on $E_\tau$. Then the monodromy operates by*

$$\widetilde{T}(e_1) = e_1 \ , \ \widetilde{T}(e_2) = e_2 - Ne_1 \ .$$

*Proof.* Consider an explicit description of the universal elliptic curve, e.g. [HuK] §7. In this description it is clear that $\widetilde{T}$ operates by left multiplication on the stalk of the local system $\mathcal{H}$ over $\mathrm{SL}_2(\mathbb{Z}) \setminus \mathfrak{H}$ with basis $1, -\tau$. This gives the above formulae. □

We want to compute the boundary morphism for the derived functor of $i'^* j'_*$ applied to the short exact sequence of sheaves

$$0 \to \mathcal{H}^\vee \otimes t^* \mathcal{L}\mathrm{og} \to t^* \mathcal{P}\mathrm{ol}^{\mathrm{coh}} \to \mathbb{Q} \to 0 \ .$$

This is nothing but the boundary morphism for the functor $H^0(\Gamma_\infty, \cdot)$ (group cohomology) applied to the same sequence considered as representations of $\Gamma_\infty$.

**Proposition C.3.2.** *Let $T = \widetilde{T}^{-1}$. Let $\widetilde{e}_i \in t^* \mathcal{P}\mathrm{ol}$ be lifts of $e_i$ under the projection $t^* \mathcal{P}\mathrm{ol} \to \mathcal{H}$.*

**a)** *We have* $T(e_1) = e_1$ *and* $T(e_2) = e_2 + Ne_1$*. The lift* $\widetilde{e}_1$ *can be assumed invariant under the monodromy.*

**b)** *There is an element* $\widetilde{u} \in t^* \mathcal{L}\mathrm{og}$ *which is invariant under the map* $[N+1]_*$ *induced by multiplication on the elliptic curve. Via the map* $\widehat{\mathfrak{U}}_{\mathcal{H}} \otimes t^* \mathcal{L}\mathrm{og} \to t^* \widehat{\mathfrak{U}}_{\mathcal{H}}$*, this induces an isomorphism* $\widehat{\mathfrak{U}}_{\mathcal{H}} \xrightarrow{\otimes \widetilde{u}} t^* \mathcal{L}\mathrm{og}$*. It agrees with* $\varrho_t$ *(see corollary A.2.6).*

**c)** *In* $t^* \mathcal{P}\mathrm{ol} / e_1 t^* \mathcal{L}\mathrm{og}$ *we have the equality*

$$T(\widetilde{e}_2) = \widetilde{e}_2 + N\widetilde{e}_1 + \sum_{j \geq 0} r_j e_2^j \widetilde{u}$$

*with*

$$r_j = \frac{N}{j!} B_{j+1} \left( \frac{a}{N} \right) .$$

Before we prove the proposition, we deduce the theorem.

*Proof.* (of theorem C.1.1) The defining diagram for $\mathcal{P}\mathrm{ol}^{\mathrm{coh}}$ was

$$\begin{array}{ccccccccc}
0 & \longrightarrow & \mathcal{H} \otimes t^* \mathcal{L}\mathrm{og} & \longrightarrow & \mathcal{P}\mathrm{ol}^{\mathrm{coh}} & \longrightarrow & \mathbb{Q} & \longrightarrow & 0 \\
 & & \downarrow{=} & & \downarrow & & \downarrow{\scriptstyle 1 \mapsto e_1^\vee \otimes e_1 + e_2^\vee \otimes e_2} & & \\
0 & \longrightarrow & \mathcal{H} \otimes t^* \mathcal{L}\mathrm{og} & \longrightarrow & \mathcal{H} \otimes t^* \mathcal{P}\mathrm{ol} & \longrightarrow & \mathcal{H}^\vee \otimes \mathcal{H} & \longrightarrow &
\end{array}$$

In order to compute the residue of the first line it is enough to compute the effect of the residue map of the second line on $e_1^\vee \otimes e_1 + e_1^\vee \otimes e_2$. We do this in the explicit complexes computing cohomology of $\Gamma_\infty$:

$$\begin{array}{ccccccccc}
0 & \longrightarrow & \mathcal{H}^\vee \otimes t^* \mathcal{L}\mathrm{og} & \longrightarrow & \mathcal{H}^\vee \otimes t^* \mathcal{P}\mathrm{ol} & \longrightarrow & \mathcal{H}^\vee \otimes \mathcal{H} & \longrightarrow & 0 \\
 & & \downarrow{\scriptstyle \mathrm{id} - T} & & \downarrow{\scriptstyle \mathrm{id} - T} & & \downarrow{\scriptstyle \mathrm{id} - T} & & \\
0 & \longrightarrow & \mathcal{H}^\vee \otimes t^* \mathcal{L}\mathrm{og} & \longrightarrow & \mathcal{H}^\vee \otimes t^* \mathcal{P}\mathrm{ol} & \longrightarrow & \mathcal{H}^\vee \otimes \mathcal{H} & \longrightarrow & 0
\end{array}$$

(Note that the standard complex would be with $\widetilde{T} - \mathrm{id} = \widetilde{T}(\mathrm{id} - T)$ but this is quasi-isomorphic to the one above via multiplication with $\widetilde{T}$, which operates trivially on invariants and coinvariants.) Note that

$$Te_1^\vee = e_1^\vee - Ne_2^\vee \ , \ Te_2^\vee = e_2^\vee \ .$$

Hence the residue of $e_1^\vee \otimes e_1 + e_2^\vee \otimes e_2$ is given by

$$\begin{aligned}(\mathrm{id} - T)(e_1^\vee \otimes \widetilde{e}_1 + e_2^\vee \otimes \widetilde{e}_2) &= e_2^\vee \otimes (-T\widetilde{e}_2 + \widetilde{e}_2 + N\widetilde{e}_1) \\ &= \sum_{j\geq 0} -r_j e_2^\vee \otimes e_2^j \widetilde{u} \\ &\quad \in (\mathcal{H}^\vee \otimes t^* \mathcal{P}\mathrm{ol})/(T - \mathrm{id})\mathcal{H}^\vee \otimes t^* \mathcal{L}\mathrm{og} \ .\end{aligned}$$

As it should this an element of

$$H^1(\Gamma_\infty, \mathcal{H}^\vee \otimes t^* \mathcal{L}\mathrm{og}) = \mathcal{H}^\vee \otimes t^* \mathcal{L}\mathrm{og} / (\mathrm{id} - T)\mathcal{H}^\vee \otimes t^* \mathcal{L}\mathrm{og} \ .$$

We view it as an element of $H^1(\Gamma_\infty, \prod_{k\geq 0} \mathrm{Sym}^k \mathcal{H})$ via the isomorphism in c) of the proposition, i.e., via $\varrho_t$. Application of the projection from lemma 2.2.3 gives the final result because $\mathrm{pr}\, e_2^{\otimes k+1} = \frac{k+1}{k+2}\, e_2^{\otimes k}$

$$\mathrm{res}(\varrho_t t^* \mathcal{P}\mathrm{ol})^{k+1} = -\frac{k+1}{k+2}\, r_{k+1} = \frac{-N}{(k+2)k!} B_k\left(\frac{a}{N}\right) \ .$$

This finishes the proof of theorem C.1.1. □

It remains to show proposition C.3.2. We need an explicit description of $\mathcal{L}\mathrm{og}$ and $\mathcal{P}\mathrm{ol}$. It will be given in terms of the Lie algebra associated to fundamental group local system on the elliptic curve.

If $\gamma$ is a path from $p$ to $p'$ and $\gamma'$ a path from $p'$ to $p''$, we denote by $\gamma'\gamma$ the homotopy class of the composed path from $p$ to $p''$.

**Definition C.3.3.** *For a point $p \in E_\tau \smallsetminus 0$, ($p = x + y\tau$, $0 \leq x, y < 1$), define elements in $\pi_1(E_\tau \smallsetminus 0, p)$*

$$\begin{aligned}\gamma_1(p) &:= [p, p+1] \ , \\ \gamma_2(p) &:= [p, p+\tau] \ , \\ \phi_0(p) &:= \gamma_1(p)\gamma_2(p)\gamma_1(p)^{-1}\gamma_2(p)^{-1} \ .\end{aligned}$$

*Let $\varepsilon$ be a point in the fundamental parallelogram of $E_\tau$ which is very close to $0$. Note that $\phi_0(\varepsilon)$ is homotopic to the small loop around $0 \in E_\tau$ oriented inverse to the standard orientation. Let $\gamma_{\varepsilon,t}$ be a straight path from $\varepsilon$ to $t$ in the fundamental parallelogram. Conjugation with $\gamma_{\varepsilon,t}$ induces an isomorphism*

$$\begin{aligned}\gamma_{\varepsilon,t*} : \pi_1(E_\tau \smallsetminus 0, \varepsilon) &\to \pi_1(E_\tau \smallsetminus 0, t) \ , \\ \alpha &\mapsto \gamma_{\varepsilon,t}\alpha\gamma_{\varepsilon,t}^{-1}\end{aligned}$$

**Remark:** Application of $T$ or of $[N+1]_*$ maps $\varepsilon$ to a different point. However, these images are contained in a small simply connected neighbourhood of $\varepsilon$. We identify fundamental groups or stalks of sheaves in $\varepsilon$ with the corresponding objects at the image via any path in this neighbourhood. We write $T(\varepsilon) = \varepsilon$ and $[N+1]_*\varepsilon = \varepsilon$. Similarly, if a point $p$ is on edge of the fundamental parallelogram ($x = 0$ or $y = 0$), then the path $\gamma_1(p)$ respectively $\gamma_2(p)$ is not in $E_\tau \smallsetminus 0$. We solve this by replacing $p$ by a point $p + \varepsilon$ very close to it but inside the fundamental parallelogram. Hence $\gamma_i(p)$ depends on choices but our results do not.

For a finitely generated group $\pi$, let $\pi^{\mathrm{alg}}$ as in [Del2] 9.8 and $\mathrm{Lie}\,\pi$ the associated pro-nilpotent Lie algebra ([Del2] 9.2 and 9.3). It is characterized by the fact that the category of unipotent representations of $\pi$ is equivalent to the category of nilpotent representations of $\mathrm{Lie}\,\pi$ ([Del2] 9.4). There is a natural group homomorphism

$$\pi \to \pi^{\mathrm{alg}} .$$

As sets $\pi^{\mathrm{alg}}$ and $\mathrm{Lie}\,\pi$ are in bijection under the maps

$$\mathrm{Lie}\,\pi \underset{\overleftarrow{\ \log\ }}{\overset{\exp}{\longrightarrow}} \pi^{\mathrm{alg}} .$$

For the explicit construction and further computation rules we refer to [Bou] Ch II, §6 and Ex. §6, in particular No 5 and Ex. 4. The Hausdorff group on $\mathrm{Lie}\,\pi$ is isomorphic to $\pi^{\mathrm{alg}}$. As usual the map $\mathrm{ad}_{\widetilde{x}}$ on the Lie algebra is given by $\mathrm{ad}_{\widetilde{x}}\,\widetilde{a} = [\widetilde{x}, \widetilde{a}]$.

**Definition C.3.4.** *For $p \in E_\tau \smallsetminus 0$, let*

$$L_p := \mathrm{Lie}\,\pi_1(E_\tau \smallsetminus 0, p) .$$

*They form a local system of Lie algebras on $E$. Its central series filtration is*

$$L_p^1 := L_p \ , \ L_p^{i+1} = [L_p^i, L_p] .$$

*In $L_p$, let*

$$\begin{aligned} \widetilde{e}_i(p) &:= \log \gamma_i(p) \ , \\ \widetilde{u}_0(p) &:= \log \phi_0(p) = [\widetilde{e}_1(p), \widetilde{e}_2(p)] \ . \end{aligned}$$

*Note that $L_p^1/L_p^2 = \mathcal{H}$. The reduction of $\widetilde{e}_i(p)$ modulo $L_p^2$ is $e_i$.*

**Proposition C.3.5 ([BeL] 1.4.3).** *There is a commutative diagram of isomorphisms of local systems on $E \smallsetminus 0$*

$$\begin{array}{ccccccccc} 0 & \longrightarrow & \mathcal{L}\mathrm{og} & \longrightarrow & \mathcal{P}\mathrm{ol} & \longrightarrow & \pi^*\mathcal{H} & \longrightarrow & 0 \\ & & \downarrow & & \downarrow & & \downarrow{=} & & \\ 0 & \longrightarrow & L^2/[L^2,L^2] & \longrightarrow & L/[L^2,L^2] & \longrightarrow & \pi^*\mathcal{H} & \longrightarrow & 0 \end{array}$$

*The standard filtration on $\mathcal{L}\mathrm{og}$ is given by $W_{-i}\,\mathcal{L}\mathrm{og} = L^{i+2}/L^{i+2} \cap [L^2, L^2]$. The standard generator $1 \in \mathrm{Gr}_0\,\mathcal{L}\mathrm{og}$ is mapped to the class of $\widetilde{u}_0(p) \in L^2/[L^2, L]$.*

Now we can start with explicit computations.

**Lemma C.3.6. a)** *In $\pi_1(E_\tau \smallsetminus 0, \varepsilon)$ we have*

$$\begin{aligned} T(\gamma_1(\varepsilon)) &= \gamma_1(\varepsilon)\ , \\ T(\gamma_2(\varepsilon)) &= \gamma_2(\varepsilon)\gamma_1(\varepsilon)^N\ , \\ T(\phi_0) &= \phi_0\ , \\ T(\gamma_{\varepsilon,t}) &= \gamma_{\varepsilon,t}\gamma_1(\varepsilon)^a\ . \end{aligned}$$

**b)** *There is a commutative diagram on fundamental groups*

$$\begin{array}{ccc} \pi_1(E_\tau \smallsetminus 0, \varepsilon) & \xrightarrow{\gamma_{\varepsilon,t*}} & \pi_1(E_\tau \smallsetminus 0, t) \\ {\scriptstyle [N+1]_*}\downarrow & & \downarrow{\scriptstyle [N+1]_*} \\ \pi_1(E_\tau \smallsetminus 0, \varepsilon) & \xrightarrow{\gamma_{\varepsilon,t*}\gamma_{\varepsilon,Nt*}} & \pi_1(E_\tau \smallsetminus 0, t) \end{array}$$

*where*

$$\gamma_{\varepsilon,Nt} := \gamma_{\varepsilon,t}^{-1}([N+1]_*\gamma_{\varepsilon,t})$$

*as homotopy class in $E_\tau \smallsetminus 0$. Moreover, $[N+1]_*\phi_0 = \phi_0$ in $\pi_1(E_\tau \smallsetminus 0, \varepsilon)$.*

*Proof.* $T$ on fundamental groups is the same as the map induced by the automorphism $T$ of $E_\tau$ which maps $\tau \mapsto \tau + N$. In the cover $\mathbb{C} \smallsetminus (\mathbb{Z} + \tau\mathbb{Z})$ of $E_\tau \smallsetminus 0$, we can write down all paths and the homotopies between them are immediate. As $\phi_0$ is a small loop it is mapped to the same homotopy class by $[N + 1]_*$. □

For $\widetilde{x} \in L_p$ we write

$$e^{\widetilde{x}} := \sum_{k \geq 0} \frac{1}{k!}\,\mathrm{ad}_{\widetilde{x}}^k$$

as maps on $L_p$. For $\gamma \in \pi_1(E_\tau \smallsetminus 0, p)$, the map $\gamma_*$ on fundamental groups induces the map $e^{\log\gamma}$ on $L_p$ ([Bou] Ch.II §6 Ex. 1).

**Proposition C.3.7.** *Let*

$$\widetilde{u} := \gamma_{\varepsilon,t*} e^{\left(-\frac{a}{N}\widetilde{e}_2(\varepsilon)-\frac{b}{N}\widetilde{e}_2(\varepsilon)\right)}\widetilde{u}_0(\varepsilon)$$

*in* $\mathcal{L}\mathrm{og}_t = L_t^2/[L_t^2, L_t^2]$. *Then*

$$[N+1]_*\widetilde{u} = \widetilde{u} \ ,$$
$$\widetilde{u} \equiv \widetilde{u}_0(t) \quad \mod W_{-1}\,\mathcal{L}\mathrm{og} = [L_t, L_t^2] \cap L_t^2 \ .$$

*Proof.* Recall that $t = \frac{a}{N}\tau + \frac{b}{N}$. As maps on $\mathcal{L}\mathrm{og}$, we have the formula

$$e^{\widetilde{e}_2(\varepsilon)}e^{\widetilde{e}_1(\varepsilon)} = e^{\widetilde{e}_1(\varepsilon)}e^{\widetilde{e}_2(\varepsilon)} = e^{\widetilde{e}_2(\varepsilon)+\widetilde{e}_1(\varepsilon)}$$

because the difference will be in $[L_\varepsilon^2, L_\varepsilon^2]$. The path $\gamma_{\varepsilon,Nt}$ defined in lemma C.3.6 is a product of $a$ copies of $\gamma_2(\varepsilon)$ and $b$ copies of $\gamma_1(\varepsilon)$ in some order. However, on $\mathcal{L}\mathrm{og}$

$$\gamma_{\varepsilon,Nt*} = e^{a\widetilde{e}_2(\varepsilon)+b\widetilde{e}_1(\varepsilon)} \ .$$

With these rules and C.3.6 b) we have

$$\begin{aligned}
[N+1]_*\widetilde{u} &= [N+1]_*\ \gamma_{\varepsilon,t*}\ e^{-\frac{a}{N}\widetilde{e}_2(\varepsilon)-\frac{b}{N}\widetilde{e}_1(\varepsilon)}\ \widetilde{u}_0(\varepsilon)\\
&= \gamma_{\varepsilon,t*}\ e^{a\widetilde{e}_2(\varepsilon)+b\widetilde{e}_1(\varepsilon)}\ [N+1]_*\ e^{-\frac{a}{N}\widetilde{e}_2(\varepsilon)-\frac{b}{N}\widetilde{e}_1(\varepsilon)}\ \widetilde{u}_0(\varepsilon)\\
&= \gamma_{\varepsilon,t*}\ e^{a\widetilde{e}_2(\varepsilon)+b\widetilde{e}_1(\varepsilon)}\ e^{-\frac{(N+1)a}{N}\widetilde{e}_2(\varepsilon)-\frac{(N+1)b}{N}\widetilde{e}_1(\varepsilon)}\ \widetilde{u}_0(\varepsilon)\\
&= \gamma_{\varepsilon,t*}\ e^{-\frac{a}{N}\widetilde{e}_2(\varepsilon)-\frac{b}{N}\widetilde{e}_1(\varepsilon)}\ \widetilde{u}_0(\varepsilon)
\end{aligned}$$

The second equality is clear from the definition. □

**Corollary C.3.8.** *Write* $e_1^i e_2^j = e_1^{\otimes i} e_2^{\otimes j} \in \mathrm{Sym}^{i+j}\,\mathcal{H}$. *Let*

$$e_1^i e_2^j \widetilde{u} := \mathrm{ad}_{\widetilde{e}_1}^i\,\mathrm{ad}_{\widetilde{e}_2}^j\,\widetilde{u} \in \mathcal{L}\mathrm{og}_t \ .$$

*This gives an isomorphism*

$$\prod_{k\geq 0}\mathrm{Sym}^k\,\mathcal{H} \cong \mathcal{L}\mathrm{og}_t$$

*of* $\mathcal{L}\mathrm{og}_0$*-modules. It is equal to the isomorphism* $\varrho_t$ *of corollary A.2.6.*

*Proof.* The space of elements in $t^* \mathcal{L}\mathrm{og}$ which are invariant under $[N+1]_*$ is one-dimensional. Hence a splitting by a $[N+1]$-invariant element is algebraic. The isomorphism $\varrho_t$ of mixed sheaves is uniquely determined as section of the weight filtration. □

Let $\overline{\mathcal{P}\mathrm{ol}}_t = \mathcal{P}\mathrm{ol}_t \,/\, \mathrm{ad}_{\widetilde{e}_1}(\mathcal{L}\mathrm{og}_t) = \mathcal{P}\mathrm{ol}_t \,/ e_1\, \mathcal{L}\mathrm{og}_t$. We write

$$z^i = \mathrm{ad}^i_{\widetilde{e}_2(t)}$$

on $\overline{\mathcal{P}\mathrm{ol}}_t$. The complete vector space $\overline{\mathrm{Pol}}_t$ has the basis elements $\widetilde{e}_2(t)$ and $z^i\widetilde{e}_1(t)$ for $i \geq 0$. All $z^i\widetilde{e}_1(t)$ commute.

**Lemma C.3.9.** *Let $V$ be in the ideal generated by $\widetilde{e}_1$ of $L_\varepsilon$ and $U \in L_\varepsilon$. Then in $\overline{\mathcal{P}\mathrm{ol}}_\varepsilon$:*

$$\log \exp(U+V) = \log\left(\exp\left(\frac{e^U - 1}{\mathrm{ad}_U} V\right)\exp U\right) .$$

*Proof.* By [Bou] Ch. II §6 Prop. 5, we have

$$(*) \qquad \log(\exp(U+V)\exp(-U)) = \sum_{n\geq 0} \frac{1}{(n+1)!}\,\mathrm{ad}^n_U V$$

up to relations in

$$\sum_{m\geq 0}\sum_{n\geq 2} A^{m,n}(\{U,V\}) \cap L(\{U,V\}) = \sum_{m\geq 0}\sum_{n\geq 2} P_{n+m}A^{m,n}(\{U,V\})$$

with notations in loc. cit. and the projection $P_{n+m}$ as in the proof of loc. cit. Thm 2. It is easy to see that

$$P_{n+m}A^{m,n}(\{U,V\}) \subset \mathrm{ad}_V[U,V] \subset \mathrm{ad}_{\widetilde{e}_1} \mathcal{L}\mathrm{og}_\varepsilon$$

for $n \geq 2$. Hence $(*)$ holds in $\overline{\mathcal{P}\mathrm{ol}}$. The formula in the lemma follows from $(*)$ by using the Hausdorff group law on $\overline{\mathcal{P}\mathrm{ol}}$. □

**Proposition C.3.10.** *In $L_t$, we have*

$$(1) \qquad \gamma_{\varepsilon,t*}(\widetilde{e}_i(\varepsilon)) = \widetilde{e}_i(t) ,$$

$$(2) \qquad T(\widetilde{e}_1(t)) = \widetilde{e}_1(t) ,$$

$$(3) \qquad T([\widetilde{e}_2(\varepsilon), \widetilde{e}_1(\varepsilon)]) = [\widetilde{e}_2(\varepsilon), \widetilde{e}_1(\varepsilon)]$$

*and in* $\overline{\mathcal{P}\mathrm{ol}}$*:*

$$\widetilde{u}_0(\varepsilon) = (1 - e^z)\, \widetilde{e}_1(\varepsilon)\ , \tag{4}$$

$$\widetilde{u} = \gamma_{\varepsilon,t*}\, e^{-\frac{a}{N}z}(1 - e^z)\, \widetilde{e}_1(\varepsilon)\ , \tag{5}$$

$$T(\widetilde{e}_2(\varepsilon)) = \widetilde{e}_2(\varepsilon) + \frac{Nze^z}{e^z - 1}\, \widetilde{e}_1(\varepsilon)\ , \tag{6}$$

$$T(\widetilde{e}_2(t)) - N\widetilde{e}_1(t) - \widetilde{e}_2(t) = N \sum_{j\geq 0} B_{j+1}\left(\frac{a}{N}\right)\frac{z^j}{j!}\, \widetilde{u} \tag{7}$$

*Proof.* In this proof, we drop $(\varepsilon)$ from the notation. By definition

$$\gamma_{\varepsilon,t*}(\log \gamma_i) = \log \gamma_{\varepsilon,t*}\gamma_i = \log \gamma_t(t)\ .$$

(2) and (3) are clear from lemma C.3.6 a). In the following computations we write $=$ for equalities in the Lie algebra and $\equiv$ for congruences modulo $\mathrm{ad}_{\widetilde{e}_1}(\mathcal{L}\mathrm{og}_t)$. Again by definition and the remark before C.3.7

$$\begin{aligned}\widetilde{u}_0 =& \log\left(\exp \widetilde{e}_1 \exp \widetilde{e}_2 \exp^{-1} \widetilde{e}_1 \exp^{-1} \widetilde{e}_1\right)\\ =& \log\left(\exp \widetilde{e}_1 \exp(-e^{\widetilde{e}_2}\widetilde{e}_1)\right)\\ =& \log\exp((1 - e^{\widetilde{e}_2})\widetilde{e}_1) = (1 - e^{\widetilde{e}_2})\widetilde{e}_1\ .\end{aligned}$$

This is (4). From the definition we have

$$\widetilde{u} \equiv \gamma_{\varepsilon,t*}\, e^{-\frac{a}{N}z}\widetilde{u}_0$$

which combined with (4) gives (5). By C.3.6 a)

$$\begin{aligned}T(\widetilde{e}_2) &= \log\left(\exp \widetilde{e}_2 \exp N\widetilde{e}_1 \exp^{-1} \widetilde{e}_2 \exp \widetilde{e}_2\right)\\ &= \log\left(\exp(Ne^{\widetilde{e}_2}\widetilde{e}_1) \exp \widetilde{e}_2\right)\\ &\equiv \log\exp\left(\widetilde{e}_2 + \frac{N\,\mathrm{ad}_{\widetilde{e}_2}\ e^{\widetilde{e}_2}}{e^{\widetilde{e}_2} - 1}\widetilde{e}_1\right)\end{aligned}$$

by lemma C.3.9. This is (6). Finally for (7),

$$\begin{aligned}
T(\widetilde{e}_2(t)) &= \log(T\gamma_2(t)) \\
&= \log\left(T(\gamma_{\varepsilon,t})T(\gamma_2)T(\gamma_{\varepsilon,t}^{-1})\right) \\
&= \log\left(\gamma_{\varepsilon,t}\gamma_1^a T(\gamma_2)\gamma_1^{-a}\gamma_{\varepsilon,t}^{-1}\right) \qquad \text{by C.3.2 a)} \\
&= \gamma_{\varepsilon,t*}\log\left(\exp(e^{a\widetilde{e}_1}T\widetilde{e}_2)\right) \\
&= \gamma_{\varepsilon,t*}\left(e^{a\widetilde{e}_1}T\widetilde{e}_2\right) \\
&\equiv \gamma_{\varepsilon,t*}(T\widetilde{e}_2 + a[\widetilde{e}_1, T\widetilde{e}_2]) \\
&= \gamma_{\varepsilon,t*}(T\widetilde{e}_2 - aT[\widetilde{e}_2, \widetilde{e}_1]) \\
&= \gamma_{\varepsilon,t*}\left(\widetilde{e}_2 + \frac{Nze^z}{e^z-1}\widetilde{e}_1 - az\widetilde{e}_1\right) \quad \text{by (6) and (3).}
\end{aligned}$$

Hence

$$\begin{aligned}
T(\widetilde{e}_2(t)) - N\widetilde{e}_1(t) - \widetilde{e}_2(t) &= \gamma_{\varepsilon,t*}\left(N\frac{ze^z}{e^z-1} - az - N\right)\widetilde{e}_1 \\
&\overset{(5)}{=} N\left(\frac{ze^z}{e^z-1} - \frac{az}{N} - 1\right)\frac{e^{\frac{a}{N}z}}{1-e^z}\widetilde{u} \\
&= N\partial_z\left(\frac{z\, e^{\frac{a}{N}z}}{e^z-1}\right)\widetilde{u}
\end{aligned}$$

And this is equal to the series in (7) by definition of Bernoulli polynomials.
□

Collecting everything, we have proved proposition C.3.2.